\documentclass[11pt]{amsart}
\usepackage{url,graphicx,tabularx,array,geometry}
\usepackage{mathabx,epsfig}
\usepackage{graphicx}
\usepackage{amssymb}
\usepackage{fancyhdr}
\usepackage{wrapfig}
\usepackage{epstopdf}
\usepackage{amsmath}
\usepackage[labelfont=bf]{caption}
\usepackage[hidelinks]{hyperref}
\usepackage{amsthm}
\usepackage{bbm}
\usepackage[T1]{fontenc}
\usepackage{yfonts}
\usepackage{harpoon}
\usepackage[all]{xy}
\usepackage{tikz}
\usepackage{enumitem}

\usepackage{rotating}
\usetikzlibrary{decorations.markings}
\tikzstyle{vertex}=[circle, draw, inner sep=0pt, minimum size=6pt]
\usepackage{hyperref}

\usepackage{subfigure}

\DeclareMathOperator{\Aut}{Aut}
\DeclareMathOperator{\Out}{Out}

\DeclareMathOperator{\Fit}{Fit}
\DeclareMathOperator{\PSL}{PSL}

\DeclareMathOperator{\GL}{GL}
\DeclareMathOperator{\SL}{SL}

\def\acts{\mathrel{\reflectbox{$\righttoleftarrow$}}}

\newcommand{\CP}{\mathcal{P}}

\newcommand{\CQ}{\mathcal{Q}}
\newcommand{\CR}{\mathcal{R}}
\newcommand{\CI}{\mathcal{I}}
\newcommand{\CD}{\mathcal{D}}
\newcommand{\CO}{\mathcal{O}}

\newcommand{\BZ}{\mathbb{Z}}

\newcommand{\BC}{\mathbb{C}}

\mathchardef\hyp="2D

\makeatletter
\def\imod#1{\allowbreak\mkern10mu({\operator@font mod}\,\,#1)}
\makeatother

\theoremstyle{plain}

\newtheorem{theoremN}{Theorem}[section]

\newtheorem{corollaryN}[theoremN]{Corollary}
\newtheorem{lemmaN}[theoremN]{Lemma}

\theoremstyle{definition}

\theoremstyle{remark}

\newtheorem*{note}{Note}
\theoremstyle{definition}

\newcommand{\tightoverset}[2]{%
  \mathop{#2}\limits^{\vbox to -.5ex{\kern-0.75ex\hbox{$#1$}\vss}}}

\newcommand{\tikzoverset}[2]{%
  \tikz[baseline=(X.base),inner sep=0pt,outer sep=0pt]{%
    \node[inner sep=0pt,outer sep=0pt] (X) {$#2$};
    \node[yshift=1pt] at (X.north) {$#1$};
}}

\newcommand{\overharpo}{\tikzoverset{\rightharpoonup}}
\newcommand{\loverharpo}{\tikzoverset{\leftharpoonup}}

\begin{document}
\title[A Classification of Prime Graphs of Pseudo-solvable Groups]{A Classification of Prime Graphs of Pseudo-solvable Groups}

\author[Huang, Keller, Kissinger, Plotnick, Roma, Yang]{Ziyu Huang, Thomas Michael Keller, Shane Kissinger, Wen Plotnick, Maya Roma, Yong Yang}

\address{Ziyu Huang, Department of Mathematics, Boston College, 140 Commonwealth Ave, Chestnut Hill, MA 02467, USA}
\email{huangaaf@bc.edu}

\address{Thomas M. Keller, Department of Mathematics, Texas State University, 601 University Drive, San Marcos, TX 78666-4616, USA}
\email{keller@txstate.edu}

\address{Shane Kissinger, Department of Mathematics, Science Center Room 325, Harvard University, 1 Oxford Street, Cambridge, MA 02138, USA}
\email{skissinger@college.harvard.edu}

\address{Wen Plotnick, Department of Mathematics, 2074 East Hall, University of Michigan Ann Arbor, 530 Church Street, Ann Arbor, MI 48109-1043, USA}
\email{plotnw@umich.edu}

\address{Maya Roma, Department of Mathematics, 970 Evans Hall 3840, University of California Berkeley, Berkeley, CA 94720-3840 USA}
\email{mayaroma@berkeley.edu}

\address{Yong Yang, Department of Mathematics, Texas State University, 601 University Drive, San Marcos, TX 78666-4616, USA}
\email{yang@txstate.edu}

\keywords {prime graph, $A_5$, $\SL(2,5)$}

\begin{abstract} The prime graph  $\Gamma(G)$ of a finite group $G$ (also known as the Gruenberg-Kegel graph) has as its vertices the prime divisors of $|G|$, and $p\text-q$ is an edge in $\Gamma(G)$ if and only if $G$ has an element of order $pq$. Since their inception in the 1970s these graphs have been studied extensively; however, completely classifying the possible prime graphs for larger families of groups remains a difficult problem.  For solvable groups such a classification was found in 2015. In this paper we go beyond solvable groups for the first time and characterize prime graphs of a more general class of groups we call pseudo-solvable. These are groups whose composition factors are either cyclic or $A_5$. The classification is based on two conditions: the vertices $\{2,3,5\}$ form a triangle in $\overline\Gamma(G)$ or $\{p,3,5\}$ form a triangle for some prime $p\neq 2$. 


\end{abstract}

\maketitle
\section{Introduction}
In this paper, we continue the study of prime graphs of finite groups. For a finite group $G$, its prime graph $\Gamma(G)$ is a simple graph with the set of vertices $\pi(G)$, the set of prime divisors of $|G|$. For any two primes $p$, $q\in\pi(G)$, $p\hyp q$ is an edge in $\Gamma(G)$ if and only if there exists an element in $G$ of order $pq.$ In 2015, Gruber, Keller, Lewis, Naughton, and Strasser studied the prime graphs of solvable groups (which are the groups with only cyclic composition factors) and obtained the following characterization:
\begin{lemmaN}\cite[Theorem 2.10]{GRUBER2015397}\label{gruber theorem}
    A simple graph is isomorphic to the prime graph of a solvable group if and only if the complement of the graph is 3-colorable and triangle-free.
\end{lemmaN}
This result was somewhat surprising in that one had not expected that there would be such a neat description for prime graphs of such a large class of groups. In this paper we take the first step
into the realm of non-solvable groups by completely characterizing the prime graphs of what we call 
pseudo-solvable groups. These are finite groups whose composition factors are either cyclic or isomorphic to $A_5$, the alternating group on five elements. (In particular, they are a generalization of solvable groups.) In a way pseudo-solvable groups are as close to solvable (without having to be solvable) as
possible: Their only possible non-cyclic composition factors are isomorphic to the smallest non-abelian simple group. So na\"{\i}vely one might think that their prime graphs must be differing as little as possible
from the prime graphs of solvable groups, but this is far from true. In forthcoming work
\cite{edwardskellerpesaklatha}, which builds
on the ideas and techniques developed here, the prime graphs of many more classes of groups are classified, where $A_5$ is replaced by other arithmetically small non-abelian simple groups, and it turns
out that the prime graphs are much less complicated than the ones for pseudo-solvable groups 
studied here. This shows that despite (or probably because of) its smallness, $A_5$ is "extra hard",
and the reason for this is easily explained. Edges in the complement of the prime graph often
come from Frobenius actions, and it is well-known that any non-solvable Frobenius complement must contain
a section of $\SL(2,5)$ which, oourse, has $A_5$ as a factor group. Thus $A_5$ is the only simple group
which in the guise of $\SL(2,5)$ can be part of a Frobenius action. This essentially is what complicates
matters in this paper and yields a richer structure of possible prime graphs than when allowing 
other "small" simple groups.\\
In studying pseudo-solvable groups, we are not starting from zero.
In recent work Florez, Higins, Huang, Keller, Shen, and Yang \cite{dawei2021prime} started the study of prime graphs of pseudo-solvable groups. They had two results. The first is a generalization of Lemma~\ref{gruber theorem}:\begin{lemmaN}\cite[Theorem 6.6]{dawei2021prime}\label{dawei theorem 1}
    If $G$ is a pseudo-solvable group and $\overline\Gamma(G)$ is the complement of its prime graph, then $\overline\Gamma(G)$, with the edge $3\hyp 5$ (if it exists) removed, is 3-colorable and triangle-free.
\end{lemmaN}
This result implies that all triangles in the graph $\overline\Gamma(G)$ have a base edge $3\hyp5$. They also gave a further description:\begin{lemmaN}\cite[Theorem 6.8]{dawei2021prime}\label{dawei theorem 2}
    If $G$ is a pseudo-solvable group, then one of the following holds for $\overline\Gamma(G)$:\begin{enumerate}
        \item The subgraph $\overline\Gamma(G)[\{2,3,5\}]$ is not a triangle;
        \item For any prime $p$ in $\pi(G)\setminus\{2,3,5\}$, $2\hyp p$ is not an edge in $\overline\Gamma(G)$. 
    \end{enumerate}
\end{lemmaN}
In this paper, we build on Lemma~\ref{dawei theorem 2} to give a complete classification of prime graphs of pseudo-solvable groups $G$. If $\overline\Gamma(G)$ is triangle-free and 3-colorable, then it can also be realized by some solvable group $G'$ by Lemma~\ref{gruber theorem} so these graphs are not of interest to us here any more, and we refer to these graphs as "trivial". Instead, we will discuss the following two categories of $\Gamma(G)$: 
\begin{enumerate}
    \item $\overline\Gamma(G)$ is 3-colorable and the vertices $\{2,3,5\}$ form a triangle in $\overline\Gamma(G)$;
    \item $\overline\Gamma(G)$ is 3-colorable and the vertices $\{p,3,5\}$ form a triangle for some prime $p\neq 2$ in $\overline\Gamma(G)$.
\end{enumerate}
    Theorem~\ref{3 colorable} shows that all the complement prime graphs of pseudo-solvable groups are 3-colorable, so the above two categories cover all the nontrivial cases of prime graphs of pseudo-solvable groups. Moreover, Theorem~\ref{exact characterization} tells us that the two categories are disjoint from each other. Theorems~\ref{235 triangle main theorem} and ~\ref{category 2 description} then give graph theoretic descriptions of the two categories, respectively. Both descriptions are necessary and sufficient. So these two theorems are the desired classification of the prime graohs of pseudo-solvable groups; unfortunately, they are a little too technical to state here in the introduction. We do, however, derive a neat consequence from our work, which appears below as Corollary~\ref{cor:complete classification} and which summarizes our graph-theoretical classification of pseudo-solvable groups. It is the following result.\\

        \begin{theoremN}\label{cor:complete classification}
        Let $G$ be a pseudo-solvable group and $\Gamma = \Gamma(G)$, then $\Gamma$ is in exactly one of the following categories.
        \begin{enumerate}
            \item $\overline\Gamma$ is 3-colorable and triangle-free;
            \item $\overline\Gamma$ is 3-colorable and the $\{2,3,5\}$-subgraph is a triangle. For an arbitrary prime divisor $r$, all the $\{r,2,3,5\}$-subgraphs are of types in Figure~\ref{fig: possible graphs of case 1} below;
            \item $\overline\Gamma$ is 3-colorable and there exists an odd prime $p$ such that the $\{p,3,5\}$-subgraph is a triangle. For an arbitrary prime divisor $r$, all the $\{r,2,3,5\}$-subgraphs are of types in Figure~\ref{fig: possible graphs of case 2} below;
            \item $\overline\Gamma$ is 3-colorable and there exists an odd prime $p$ such that the $\{p,3,5\}$-subgraph is a triangle. For an arbitrary prime divisor $r$, all the $\{r,2,3,5\}$-subgraphs are of types in Figure~\ref{fig: conditioned graphs of case 2} below.
        \end{enumerate}
    \end{theoremN}

    The proof of our classification results 
    Theorems~\ref{235 triangle main theorem} and ~\ref{category 2 description}
    both will have two parts. We will refer to deriving properties of the graph $\Gamma(G)$, which it needs to satisfy, as the "forward direction", and we will refer to constructing a pseudo-solvable group $G$ based on a given graph $\Gamma$ such that $\Gamma = \Gamma(G)$ as the "backward direction."\\

    We close this introduction by giving a brief outline of the paper and establishing some notation. In Section 2, we will introduce some preliminary results. In Section 3 and Section 4, we will classify the prime graphs in Category~(1) and Category~(2), respectively. Here are some conventions and notations that will be used. All groups are finite and all graphs are simple. The complement of a graph $\Gamma$ is denoted $\overline\Gamma$. Given a subset of vertices $\pi\subset V(\Gamma)$, $\Gamma[\pi]$ is the induced subgraph restricting to the vertices in $\pi$ and $\Gamma[\pi']$ is the induced subgraph restricting to the vertices that are not in $\pi$. $\pi'$ denotes the subset of vertices $V(\Gamma)\setminus\pi$. The prime graph of a group $G$ is denoted $\Gamma(G)$. When the context is clear, when we say that $2\hyp p$ is an edge we automatically assume $p\in\pi(G)\setminus\{2,3,5\}.$ 
    A pseudo-solvable group is \emph{strictly pseudo-solvable} if it is not solvable. If $G$ is an extension of $K$ by $N$, we write $G = N.K$, i.e there exists a normal subgroup $N\triangleleft G$ such that $G/N\simeq K$. If the extension splits, we write $G = N\rtimes K$, i.e $G$ is a semidirect product of $K$ and $N$ equipped with a group action $K\acts N$. Unless stated otherwise, we denote the group identity $1$. We also denote the trivial group $1.$ For any other unspecified notations, we inherit them from \cite{dawei2021prime} and \cite{GRUBER2015397}. A normal series of $G$ is a series of normal subgroups $1\triangleleft N_1\triangleleft\cdots\triangleleft N_m = G$ such that $N_i\triangleleft G$ for any $i$. We refer the quotients of consecutive terms as the \emph{factors} or \emph{sections} of $G$. We call the last factor $N_m/N_{m-1}$ the \emph{top factor} of the normal series. For each $N_i$, we refer all terms $N_j, j<i$, as the \emph{previous terms} and all terms $N_k, k>i$, as the \emph{later terms.} We abbreviate the words "without loss of generality" as WLOG. \\

    
\maketitle
\section{Preliminary results}
    We briefly introduce some more notation and present main results in \cite{dawei2021prime} and \cite{GRUBER2015397}. We first make the following definitions about Frobenius groups. 
        A group action $A\acts N$ is \emph{fixed point free} if for any $a\in A, a\neq 0$, $a$ has no nontrivial fixed points acting on $N$. Such an action is also called \emph{Frobenius}. 
        A group $G$ is \emph{Frobenius} if $G = N\rtimes K$ such that the action $K\acts N$ is Frobenius. $N$ is called a \emph{Frobenius kernel} of $G$ and $K$ is called a \emph{Frobenius complement} of $G$. 
        If $G$ is a \emph{Frobenius group} with the Frobenius action $K\acts N$, and if $K$ is a $q$-group and $N$ is a $p$-group, then we say $G$ is a \emph{Frobenius group of type $(p,q)$}.
        Let $F_1=F_1(G):= \Fit(G)$ be the Fitting subgroup of $G$ and $F_2=F_2(G)$ be the second Fitting subgroup such that $F_2/F_1 = \Fit(G/F_1)$. A solvable group $G$ is 2-Frobenius if $F_2$ is a Frobenius group with action $F_2/F_1\acts F_1$ and $G/F_1$ is a Frobenius group with action $G/F_2\acts F_2/F_1$. $F_2/F_1$ is called the \emph{upper kernel of $G$} and $F_1$ is called the \emph{lower kernel of $G$}. 
        If $G$ is a 2-Frobenius group with $F_1$ and $F_2$ specified as in the previous definition, and if $G/F_2$ and $F_1$ are $p$-groups and $F_2/F_1$ is a $q$-group, then we say $G$ is a \emph{2-Frobenius group of type $(p,q,p)$}.
        If $G = PQR$ where $P, Q, R$ are $p\hyp, q\hyp, r\hyp$groups respectively, $PQ$ is a Frobenius group of type $(p, q)$, $QR$ is either a Frobenius group of type $(q,r)$ or a 2-Frobenius group of type $(r,q,r)$, then we say $G$ is a 2-Frobenius group of type $(p,q,r)$.
    Some properties directly follow from the above definitions. If $G = N.K$ is Frobenius, then $(|N|, |K|) = 1$. The Frobenius kernel $N$ is unique and the Frobenius complement $K$ is unique up to conjugation. $N = \Fit(G)$ is a nilpotent group. The Sylow $p$-subgroup of $K$ is either cyclic or generalized quaternion. All of these results are obtained from \cite[Chapter 6]{isaacs2008finite}. Also, many results about Frobenius groups and prime graphs are described in \cite[Section 2]{GRUBER2015397} in detail.\\
    
    Let $\pi\subset\pi(G)$. A Hall $\pi$-subgroup of $G$ is a $\pi$-subgroup $H\leq G$ such that $(H, |G|/|H|) = 1$, i.e., $H$ has order of the highest powers of all prime divisors in $\pi$. We denote a Hall $\pi$-subgroup by $H_\pi$. If $\pi = \{p,q\}$ or $\pi = \{p,q,r\}$, we denote $H_\pi$ by $H_{pq}$ or $H_{pqr}$ respectively. By a theorem due to Philip Hall~\cite{hall1928note}, a finite group $G$ is solvable if and only if for any $\pi\subset \pi(G)$, $G$ has a Hall $\pi$-subgroup. Thus, when $G$ is solvable, for any $\pi\subset\pi(G)$,  $\Gamma(G)[\pi]$ is isomorphic to $\Gamma(H_\pi)$. By \cite[Theorem A]{WILLIAMS1981487}, if $\Gamma(G)$ is disconnected, then $\Gamma(G)$ has two components and $G$ is either Frobenius or 2-Frobenius. If $H$ is the Frobenius kernel (or upper Frobenius kernel) of $G$, then $\Gamma(H)$ is one of the components.\\
    
    We now define an orientation for the complement of a prime graph. First, if $G$ is solvable, Gruber et al. gave the following orientation to $\overline\Gamma(G)$ \cite[Theorem 2.8]{GRUBER2015397}. For any prime $p,q\in\pi(G)$, if $p\hyp q$ is an edge in $\overline\Gamma(G)$, then $H_{pq}$ is either a Frobenius or 2-Frobenius group. If $H_{pq}$ is Frobenius of type $(p,q)$ or 2-Frobenius of type $(p,q,p)$, we direct the edge $p\hyp q$ as $q\to p$. $\overline\Gamma(G)$ equipped with this orientation is a directed graph, and we call it the \emph{Frobenius digraph of $G$}, denoted $\overharpo\Gamma(G)$. \\
    
    
    The proof of Lemma~\ref{gruber theorem} depends on the Frobenius digraph. In the Frobenius digraph $\overharpo\Gamma(G)$, we label the vertices with in-degree 0 and nonzero out-degree by $\CO$, the vertices with nonzero in- and out-degrees by $\mathcal D$, and the remaining vertices by $\CI$. Since $\overline\Gamma(G)$ has no directed $3$-path, the labeling $\CO, \mathcal D,$ and $\CI$ gives a 3-coloring for $\overline\Gamma(G)$. We will refer to this coloring as the \emph{3-coloring for solvable graphs.} On the other hand, to show that a 3-colorable and triangle-free graph $\Gamma$ is isomorphic to $\overline\Gamma(G)$, for some pseudo-solvable group $G$, involves an explicit construction of a pseudo-solvable group $G$. That construction will be explained later in this paper and will be used for our purposes. We then mention two lemmas in \cite{dawei2021prime} which we will use later. 
    \begin{lemmaN}\cite[Lemma 6.2, Lemma 6.3]{dawei2021prime}\label{N.(A_5 timesH) subgroup}
        Let $G$ be a strict pseudo-solvable subgroup, i.e a pseudo-solvable subgroup that is not solvable, then $G$ admits a subgroup $K\simeq N.(A_5\times H)$ where $N, H$ are solvable and $(|H|, 30) = 1$. Furthermore, $|K|$ and $|G|$ have the same set of prime divisors. Thus, $\overline\Gamma(G)$ is obtained from $\overline\Gamma(K)$ by removing edges. 
    \end{lemmaN}
    The above lemma implies when we investigate graph properties that are preserved under removing edges, such as chromatic number or not containing $k$-cycles, it suffices to check them for the subgroup $K\simeq N.(A_5\times H)$ of $G.$ We also have the following lemma about normal series. \begin{lemmaN}\cite[Lemma 6.4]{dawei2021prime}\label{dawei chief series intersection}
       Let $\pi\subset\pi(G)$ and $H$ is a Hall $\pi$-subgroup of $G$. Suppose $G$ has a normal series $1 = N_0\triangleleft N_1\triangleleft N_2\triangleleft\cdots\triangleleft N_m = G$, then $H$ has an induced normal series $1 = N_0\cap H\triangleleft N_1\cap H\triangleleft N_2\cap H\triangleleft\cdots\triangleleft N_m\cap H = H$ such that $N_i\cap H$ is a Hall $\pi$-subgroup of $N_i$, for any $i\in\{0,\ldots, m\}$.
    \end{lemmaN}
    
    Florez et al. also defined orientations for edges of complements of prime graphs of pseudo-solvable groups \cite[Theorem 6.6]{dawei2021prime}. Let $G$ be a pseudo-solvable subgroup. If $G$ is solvable, then $\overline\Gamma(G)$ admits the orientation of the Frobenius digraph $\overharpo\Gamma(G)$. If $G$ is strictly pseudo-solvable, take the subgroup of $K\simeq N.(A_5\times H)$ of $G$. $K$ has subgroups $K_1\simeq N.(D_{10}\times H)$ and $K_2\simeq N.(A_4\times H)$, which both are solvable. $\overline\Gamma(K_1)$ and $\overline\Gamma(K_2)$ admit orientations. By the arguments of \cite[Theorem 6.6]{dawei2021prime}, we can assign an orientation to $\overline\Gamma(K)$ based on the orientations of $\overharpo\Gamma(K_1)$ and $\overharpo\Gamma(K_2)$. For a pair $p \text- q$ in $\overline\Gamma(K)$, if $p,q\notin\{2,3,5\}$, the direction is assigned by the rule of Frobenius digraphs in \cite{GRUBER2015397}. If $p\in \{2,3,5\}$ and $q\notin\{2,3,5\}$, then $p\to q$ is assigned. Directions $3\to 2, 2\to 5$ are assigned when the edges exists. We leave the edge $3\text-5$ with unassigned direction. Therefore, we have defined an orientation for $\overline\Gamma(K)\setminus\{3\text -5\}$. We denote it by $\overharpo\Gamma(K)$ and call it the Frobenius digraph of the group $K$. If we reverse the direction $2\to 5$ to $5\to 2$, then the graph $\overline\Gamma(K)\setminus\{3\text -5\}$ has an orientation with no directed 3-path. We denote this digraph $\loverharpo\Gamma(K)$. \\
    
    \section{Category~(1):  \{2,3,5\} forms a triangle }
    
\begin{figure}[h!]
    \centering
        \subfigure[$t_1, t_2, t_3, t_5$]{
            \begin{tikzpicture}[scale = 1]
        \node (2) at (-1,-1) {$2$};
        \node (3) at (-1,1) {$3$};
        \node (5) at (1,1) {$5$};
        \node (p) at (1,-1) {$r$};
        \path [-] (3) edge node {} (5);
        \path [-] (3) edge node {} (2);
        \path [-] (2) edge node {} (5);
            \end{tikzpicture} 
            }\;\;\;\;\;\;
        \subfigure[$t_4$]{
            \begin{tikzpicture}[scale = 1]
        \node (2) at (-1,-1) {$2$};
        \node (3) at (-1,1) {$3$};
        \node (5) at (1,1) {$5$};
        \node (p) at (1,-1) {$r$};
        \path [-] (3) edge node {} (5);
        \path [-] (3) edge node {} (2);
        \path [-] (2) edge node {} (5);
        \path [-] (5) edge node {} (p);
            \end{tikzpicture} 
            }
    \caption{Subgraphs constructed by the representations of $A_5$}
    \label{fig: possible graphs of case 1}
    \end{figure}
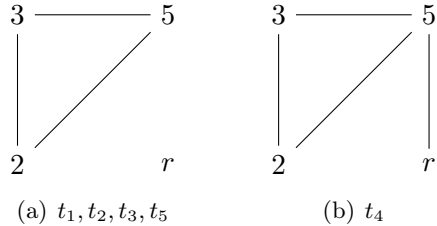

    In this section, we discuss the case when $\{2,3,5\}$ forms a triangle in the prime graph of a pseudo-solvable group $G$. By Lemma~\ref{gruber theorem}, $G$ is strictly pseudo-solvable so it must contain a section of $A_5$. We will prove that the subgraph $\overline\Gamma(G)[\{2,3,5\}']$ is still triangle-free and 3-colorable.  Moreover, for any prime $r\in\pi(G)\setminus\{2,3,5\}$, $r\text-2$ and $r\text-3$ are not edges in $\overline\Gamma(G)$. Hence the only triangle in $\overline\Gamma(G)$ is given by the vertices $\{2,3,5\}$.We first note the following lemma. 
    \begin{lemmaN}\label{plms lemma}\cite[Proposition 4.2]{plms/s3-26.4.653}
        Given $G$ a finite group and $H\triangleleft G$ such that $G/H\simeq \PSL(2,2^\alpha)$ for some $\alpha\geq 2$. If $H$ is not a $2$-group, then for any $3$-element $t\in G$, $t$ has a fixed point acting on $H$.
    \end{lemmaN}
    The case we need is when $\alpha = 2$, $\PSL(2,4)\simeq A_5$. That is, when $H$ is a $p$-group and the quotient $G/H$ is $A_5$, $t$ must fix an element of order $p$ and so eliminating the edge $3\text - p$ in $\overline\Gamma(G)$.
    \begin{lemmaN}\label{3 not adjacent}
        Let $G$ be a strictly pseudo-solvable group. If $\{2,3,5\}$ forms a triangle in $\overline\Gamma(G)$, then $\overline\Gamma(G)$ does not contain an edge $3\text- p$, for any prime $p\neq 2,5.$
    \end{lemmaN}
    \begin{proof}
        Since $G$ is strictly pseudo-solvable, it has a subgroup $K\simeq N.(A_5\times H)$ as described above. Since $\overline\Gamma(G)$ is obtained by $\overline\Gamma(K)$ removing edges, we know $\{2,3,5\}$ forms a triangle in $\overline\Gamma(K)$. It suffices to show $\overline\Gamma(K)$ does not contain an edge $3\text- p$, for any prime $p\neq 2,5.$ We claim that $(|N|, 15) = 1$. If $3\mid |N|$, we look at the subgroup $K_1\simeq N.(D_{10}\times H)$ and find that $\pi(K) = \pi(K_1)$. Thus, $\overline\Gamma(K)$ is obtained by removing edges from $\overline\Gamma(K_1)$. Since $K_1$ is solvable, $\overline\Gamma(K_1)$ does not contain any triangles, which is a contradiction. Similarly, we know $5\mid |N|$ by looking at the subgroup $K_2\simeq N.(A_4\times H)$. Thus, the claim is proved. \\
        
        For any prime $q\in \pi(H)$, there exists an element of order $3q$ so it suffices to show the statement for $N.A_5$. So without loss $K=N.A_5$ for the remainder of the proof. We choose a chief series of $K$, $1 = N_0\triangleleft N_1\triangleleft \cdots\triangleleft N_I = N \triangleleft K = N.A_5$ where $A_5$ is the top factor. If $2\nmid |N|$, then $(|N|, |A_5|) = 1$. By Schur-Zassenhaus Theorem, the extension $N.A_5$ splits so $K\simeq N\rtimes A_5$. Since $N$ is solvable, its chief factors are elementary abelian groups. Thus, for any $p\in\pi(K)\setminus \{2,3,5\}$, there exists $\alpha \in \{0,\ldots, I-1\}$ such that $N_{\alpha+1}/N_\alpha$ is an elementary abelian $p$-group. Consider $K/N_\alpha$, which is isomorphic to $(N/N_\alpha)\rtimes A_5$. $N_{\alpha+1}/N_\alpha$ is a normal subgroup of $N/N_\alpha$ so $A_5\acts N_{\alpha+1}/N_\alpha$ is a subaction of $A_5\acts N/N_\alpha$, which implies $(N_{\alpha+1}/N_\alpha)\rtimes A_5$ is a subgroup of $K/N_\alpha$. By Lemma~\ref{plms lemma}, there exists $t\in K/N_\alpha$ of order 3 such that it has a fixed point when acting on $N_{\alpha+1}/N_\alpha$. Thus, $K/N_\alpha$ has an element of order $3p$ and so does $K$. Therefore, $\overline\Gamma(K)$ does not have an edge $3\text -p$. \\
        
        We then discuss the possibility that $2\mid |N|$. Since $15\nmid |N|$, any normal series of $K\simeq N.A_5$ contains exactly one copy of $A_5$ in the factors. We first consider a special case $K = (Q\rtimes S_2).A_5$ where $S_2$ is a 2-group and $Q$ is an elementary abelian $q$-group for some other prime $q\neq 2,3,5.$ Then $K$ has a normal series $$1\triangleleft Q\triangleleft Q\rtimes S_2\triangleleft (Q\rtimes S_2).A_5.$$ Since $(|Q|, |K/Q|) = 1$, $Q$ has a complement $L\simeq S_2.A_5$. If the group action $L\acts Q$ is faithful, then $S_2$ also acts on $Q$ faithfully and thus nontrivially. For an element of order $3,$ if $\langle t\rangle $ acts on $Q\rtimes S_2$ frobeniusly, then the Frobenius kernel $Q\rtimes S_2$ is nilpotent, contradicting $S_2\acts Q$ is nontrivial. Therefore, the action $\langle t\rangle \acts Q\rtimes S_2$ is non-Frobenius, which means an element in $Q\rtimes S_2$ is fixed when conjugating by $t$. That element cannot have order $2$ since $2\text- 3$ is an edge in $\overline\Gamma(K)$. Thus, $3\text - q$ is not an edge in $\overline\Gamma(K)$. \\
        
        On the other hand, if the action $L\acts Q$ is nonfaithful, then let $1<A\triangleleft L$ be its kernel. If there is $t\in A$ of order $3$, then an element in $Q$ is fixed when conjugating by $t$ and so $3\text- q$ is not an edge in $\overline\Gamma(K)$ and we are done. Otherwise, consider $t\in L\setminus A$ of order 3. Note that $S_2A/S_2$ is a normal subgroup of $L/S_2\simeq A_5$, and so $S_2A=L$ or $S_2A=S_2$.
        If $S_2A=L$, then $t\in S_2A$. By the order of $t$, we find $t\in A$, which contradicts $t\in L\setminus A$. Thus $S_2A=S_2$, so $A\leq S_2$.
        Hence we have a faitfhul action of $L/A$ on $Q$, and $L/A$ has $A_5$ as a top composition factor, so just as in the faithful case above we
        obtain that there exists an element of order 3, $tA\in L/A$, that fixes a nontrivial element in $x\in Q$. Thus, $xt\in K$ has order divisible by $3q$ so $3\text -q$ is not an edge in $\overline\Gamma(Q)$ as desired. \\
        
        Consider the general case of $2\mid |N|$ and $15\nmid |N|$. We reduce the general case to the special case using Hall subgroups. For any $p\in \pi(K)\setminus \{2,3,5\}$, let $\pi = \{2,3,5,p\}$. Since the composition factors of $K$ are either cyclic or isomorphic to $A_5$, $K$ is $\pi$-separable so $K$ has a Hall $\pi$-subgroup, denoted $E_\pi$. By Hall's theorem, we identify $\overline\Gamma(E_\pi)$ with $\overline\Gamma(K)[\{2,3,5,p\}]$. By Lemma~\ref{dawei chief series intersection}, any chief series $1 = N_0\triangleleft N_1\triangleleft \cdots\triangleleft N_I = N.H \triangleleft K = N.(A_5\times H)$ induces a normal series by intersecting with $E_\pi$. By removing the repetitive terms, we obtain $1 = M_0\triangleleft M_1\triangleleft\cdots\triangleleft M_J = E_\pi$. For all indices $j\in\{0,\ldots, J-2\}$, $M_{j+1}/M_j$ is either an elementary abelian 2-group or an elementary abelian $p$-group, and $M_J/M_{J-1}\simeq A_5$. Pick the maximal index $r$ such that $M_{r+1}/M_r$ is an elementary abelian $p$-group. Consider $E_0 = E_\pi/M_r$. If $r = J-2$, then $E_0/(M_{J-1}/M_{J-2})\simeq A_5$ so by Lemma~\ref{plms lemma}, $E_0$ has an element of order $3p$ and so does $E_\pi$. Otherwise, $r<J-1$ so for all indices $r<j<J-1$, $M_{j+1}/M_j$ is an elementary abelian 2-group. Thus, we know $E_0 = (P\rtimes S_2).A_5$ where $P\simeq M_{r+1}/M_r$ is an elementary abelian $p$-group, $S_2$ is a Sylow 2-subgroup. By the above special case, $E_0$ has an element of order $3p$ and so does $E_\pi$. Therefore, $\overline\Gamma(K)$ does not have the edge $3\text-p$. Thus the proof is complete.
    \end{proof}
    Lemma~\ref{dawei theorem 2} and Lemma~\ref{3 not adjacent} together imply the following result.
    \begin{lemmaN}\label{2 3 not adjacent}
        Let $G$ be a strictly pseudo-solvable group. Assume $\{2,3,5\}$ forms a triangle in $\overline\Gamma(G)$. Let $r$ be a prime that divides $|G|$ but does not equal $2, 3, $ or 5, then the $\{r,2,3,5\}$-subgraph of $\overline\Gamma(G)$ must be one of the graphs described in Figure~\ref{fig: possible graphs of case 1}. 
    \end{lemmaN}
    Now we move to the graph theory aspect. We will find that the graphs in Figure~\ref{fig: possible graphs of case 1} are not only possible, but we can realize them through the irreducible representations of $A_5$. The following lemma tells us how to construct a finite group action using an ordinary representation while preserving fixed points of the representation. This is well-known, but a good reference seems to be hard to find, so we
    present the argument here for completeness.
    
    \begin{lemmaN}\label{modular representation lemma}
        Let $G$ be a finite group and $p$ prime such that $p\nmid |G|$. If $\pi: G\longrightarrow\GL(V)$ is an ordinary representation over $\BC$, then there exists a finite field $K$ of characteristic $p$ and a $K$-representation $\sigma: G\longrightarrow \GL(W)$ such that for any element $g\in G$, $\pi(g)$ and $\sigma(g)$ have the same number of fixed points. 
    \end{lemmaN}
    
    \begin{proof}
        The proof uses ideas from \cite[Section 15]{isaacs1976character} and \cite[Section 9]{aschbacher_2000}. \\
        Let $R$ be the ring of algebraic integers in $\BC$ and $M$ be a maximal ideal of $R$ containing $pR$. Thus, we have a quotient $*: R\longrightarrow R/M = F$. Let $U = \{\epsilon \in \BC: \epsilon^m = 1 \mbox{ for some integer $m$ not divisible by $p$}\}$. By \cite[Lemma 15.1]{isaacs1976character}, we know $U\subset R$ and $*|_U$ is an isomorphism from $U$ to $F^\times$; also, $F$ is an algebraically closed field of characteristic $p$. Let $\tilde R$ be the localization of $R$ at $M$, and we extend $*$ to a map $\tilde R\longrightarrow F$ canonically. \\
        
        By \cite[Theorem 15.8]{isaacs1976character}, WLOG we assume all entries of $\pi(g)$, for any $g\in G$, lie in $\tilde R$. We consider the image of $\pi$ under $*$, which is a $F$-representation $\pi^*: G\longrightarrow \GL(V^*)$ where for any $g\in G$, $\pi^*(g) = \pi(g)^*$ by applying $*$ entry by entry. Notice that since $G$ is finite, the eigenvalues of $\pi(g)$ are roots of 1 so the eigenvalues of $\pi(g)$ is contained in $U$. Let $char (g, \pi)$ be the characteristic polynomial of $\pi(g)$ and so $char(g, \pi^*)$ is the image of $char (g, \pi)$ under $*$. Since $*: U\longrightarrow F^\times$ is an isomorphism, the multiplicity of $1$ as a root in $char(g, \pi)$ is the same as the multplicity in $char(g, \pi^*)$. Thus, $\pi(g)$ has the same number of fixed points as $\pi^*(g)$. \\
        
        Let $V^* = \bigoplus_{i=1}^m Fe_i$. Since $F$ is a direct limit of finite fields of characteristic $p$, we can pick a sufficiently large finite subfield $K\subset F$ such that the image $\pi^*(G)$ lies in $\GL_m(K)\leq \GL_m(F)$ and also $W = \bigoplus_{i=1}^m Ke_i$ contains the fixed points of $\pi^*(g)$ for any $g$. Therefore, the $F$-representation $\pi^*: G\longrightarrow \GL(V^*)$ is reduced to a $K$-representation $\sigma:=\pi^*\mid_K: G\longrightarrow \GL(W)$ where all the fixed points are preserved.
    \end{proof}
    \begin{table}[h]
        \centering
        \begin{tabular}{|| c | c | c c c c c ||}
            \hline
            Index & Degrees & \multicolumn{5}{c||}{Orders} \\
            \hline
            & & 1 & 2 & 3 & 5 & 5 \\
            \hline
            $t_1$ & 1 & 0 & 0 & 0 & 0 & 0\\
            \hline
            $t_2$ & 3 & 0 & 0 & 0 & 0 & 0 \\
            \hline
            $t_3$ & 3 & 0 & 0 & 0 & 0 & 0 \\
            \hline
            $t_4$ & 4 & 0 & 0 & 0 & - & - \\
            \hline
            $t_5$ & 5 & 0 & 0 & 0 & 0 & 0 \\
            \hline
        \end{tabular}
        \caption{The values of $\det(t(g)-I)$ for irreducible representations and conjugacy class representatives of $A_5$}
        \label{tab:rep of a5}
    \end{table}
    We notice that whether a coprime action $G\acts W$ has a fixed point can be important to prime graphs since for any $q\in \pi(G)$, $q\text-p$ is an edge in $\overline\Gamma(W.G)$ if and only if the action is fixed point free. Therefore, we use GAP~\cite{gap} to check the fixed points of the irreducible representations of $A_5$. The first row of Table~\ref{tab:rep of a5} displays the orders of the representatives of conjugacy classes in $A_5$. $t_1$ through $t_5$ denotes the irreducible representations of $A_5$. The rest of the entries are calculations of $\det(t(g)-I)$ where $g$ is the conjugacy class representative of the column, $t$ is the representation of the row, and $I$ is the identity matrix of the appropriate dimension. If $\det(t(g)-I)$ equals $0$, then the entry displays "0". Otherwise the entry dislays "-". The entry values "0" or "-" signals whether the elements of that conjucacy class has fixed points through the representation. "0" implies the conjucacy class has a fixed point and "-" implies the class is fixed point free.\\
    
    In Table~\ref{tab:rep of a5}, we observe that $t_1, t_2, t_3$, and $t_5$ has fixed points for conjugacy classes of all orders. Therefore, the finite group actions that they generate through Lemma~\ref{modular representation lemma} gives complement prime graphs that are isomorphic to Figure~\ref{fig: possible graphs of case 1}(a). In $t_4$, all the elements of order 5 in $A_5$ does not have a fixed point so its finite group action gives a graph that is isomorphic to Figure~\ref{fig: possible graphs of case 1}(b), where $5\text-r$ is an edge. \\ 
    
    The next step is to construct a bigger group that contains $A_5$ as a section while its complement prime graph contains Figure~\ref{fig: possible graphs of case 1} as a subgraph. We borrow the idea from the proof of \cite[Theorem 2.8]{GRUBER2015397}. In particular, we prove an enhanced version of \cite[Lemma 1.8]{Keller1994SolvableGW}.
    \begin{lemmaN}\label{representation lemma}
        Given two finite groups $A, B$ such that $|B| = q_1\cdots q_mp_1\cdots p_n$ where$q_i, p_j$ are $m+n$ mutually distinct primes. Assume $|\Fit(B)| = q_1\cdots q_m$. Given a $\BC$-representation $\pi: A\longrightarrow\GL(U)$ and a prime $r$ such that $|A\times B|\mid r-1$, then there exists a finite field $F$ of characteristic $r$ and a $F$-representation $\rho: A\times B\longrightarrow \GL(V)$ such that the following holds:
        \begin{enumerate}
            \item For any element $a\in A$, $\pi(a)$ is fixed point free if and only if $\rho(a)$ is fixed point free;
            \item $\Fit(B)$ acts fixed point freely on $V$.
        \end{enumerate}
    \end{lemmaN}
    \begin{proof}
        Denote $G = A\times B$, $H = A\times \Fit(B)\triangleleft G$. By the order, we know $\Fit(B) = C_{q_1}\times\cdots\times C_{q_m}.$ Since $r\nmid |A|$, by Lemma~\ref{modular representation lemma}, we know there exists a finite field $F$ of characteristic $r$ and a $F$- representation $\sigma: A\longrightarrow \GL(V)$ such that for any $ a\in A$, $\pi(a)$ and $\sigma(a)$ have the same number of fixed points.  \\
        
        For any index $i\in \{1,\ldots, m\}$, since $q_i\mid r-1$, we pick a nontrivial $1$-dimensional $F$-representation $\rho_i: C_{q_i}\longrightarrow \GL(F)$, where a generator $c\in C_{q_i}$ acts on $F$ by multiplication with a $q_i$th root of unity in $F$. Consider the $F$-representation $\tilde\sigma = \sigma\otimes \rho_1\otimes\cdots\otimes \rho_n: H\longrightarrow \GL(V)$. For any $a\in A$, the number of fixed points of $a$ is preserved; for any $i, c\in C_{q_i}$, $c$ acts by scaling a nontrivial element so is fixed point free. Consider the induced representation $W= \text{Ind}_H^G V$ where we denote the $F$-representation $\rho: G\longrightarrow \GL(W)$. Then, $W$ has the direct sum decomposition $$
        W = \bigoplus_{Hy\in G/H}W_y = \bigoplus_{Hy\in G/H}V^y,
        $$
        where each $y$ is a representative in the coset $Hy$. Each $W_y$ is a $F$-subspace of $W$ that is preserved by the action of $H$, and $G$ acts on them by $g\cdot W_y = W_{yg}$. Each $V^y$ is the conjugate representation of $V$ by $y$, and since $H$ is normal, the action is by $\tilde\sigma^y(h) = \tilde\sigma(h^{y^{-1}})$. We have a canonical $H$-bijection $W_y\longleftrightarrow V^y$ so for any $ h\in H$, we know whether $h$ is fixed point free acting on $W$ by counting its fixed points acting on each $V^y$. For any $a\in A\times 1$, $\tilde\sigma^y(a)$ is fixed point free if and only if $\tilde\sigma(a)$ is fixed point free. Thus, $\rho(a)$ is fixed point free if and only if $\tilde\sigma(a)$ is fixed point free if and only if $\sigma(a)$ is fixed point free. Also, for any $b\in 1\times \Fit(B)$, since $\tilde\sigma(b)$ is fixed point free, each $\tilde\sigma^y(b)$ is fixed point free and so $\rho(b)$ is fixed point free. Therefore, $\rho: G\longrightarrow \GL(W)$ over $F$ is the desired representation.
        \end{proof}
    Notice that in Lemma~\ref{representation lemma}, if $(|A|, |B|) = 1$, given a prime $p\mid|A|$, if all $p$-elements in $A$ acts fixed point freely on $U$, then all $p$-elements in $A\times B$ act fixed point freely on $V$. This result is given by the fact that all $p$-elements in $A\times B$ lie in $A\times 1$. It guarantees that new fixed points are not generated in $\rho$ when it comes from $\pi$. In all the following cases that we apply Lemma~\ref{representation lemma}, we will know $(|A|, |B|) = 1$ and so we will assume this result implicitly. \\
    
    We now give a description for prime graphs of solvable groups in category~(1). 
    \begin{theoremN}\label{235 triangle main theorem}
        The following are equivalent for a simple graph $\Gamma$:\begin{enumerate}
            \item $\Gamma \simeq \Gamma(G)$ for some pseudo-solvable group $G$ such that $\{2,3,5\}$ forms a triangle in $\overline\Gamma(G)$;
            \item $\Gamma$ with a choice of vertices $\{a,b,c\}\in V(\Gamma)$ satisfies the following conditions: 
            \begin{enumerate}
                \item $\overline\Gamma$ is 3-colorable;
                \item $\{a,b,c\}$ forms the only triangle in $\overline\Gamma$; 
                \item For any element $r\in V(\overline\Gamma)\setminus \{a,b,c\}$, the $\{r,a,b,c\}$-subgraph of $\overline\Gamma$ is isomorphic to one of the subgraphs in Figure~\ref{fig: possible graphs of case 1} where $a$ is mapped to $2$, $b$ is mapped to $3$, and $c$ is mapped to $5$, $r$ is mapped to $r$;
                \item $\overline\Gamma$ has a choice of coloring $\CO, \CD, \CI$ such that the vertices adjacent to $c$, except for $a$ and $b$, are of one color. 
            \end{enumerate}
        \end{enumerate}
    \end{theoremN}

    \begin{proof}
        Assume condition (1), we want to prove condition (2). By Lemma~\ref{dawei theorem 1}, since $3\text-p$ is not an edge for any $p\neq 2,5$, $\overline\Gamma(G)$ has only one triangle so 2(b) holds. We then take the subgroup $K\simeq N.(H\times A_5)$ of $G$. It has an orientation $\loverharpo\Gamma(K)$ with no directed 3-paths. If we direct the edge $3\text - 5$ by $5\to 3$ and then give the new directed graph the 3-coloring $\CO, \mathcal D$, and $\CI$, we find $\overline\Gamma(K)$ is 3-colorable with $2\in\CI, 5\in\CO$, and $3\in\mathcal D$. Thus, we have showed 2(a). By Lemma~\ref{2 3 not adjacent}, 2(c) holds. We consider the subgroup $K_1\simeq N.(D_{10}\times H)$ of $K$. $K_1$ is solvable so its Frobenius digraph $\overharpo\Gamma(K_1)$ has no directed 3-paths. Due to the edge $2\to 5$ the set of vertices with edges outpointing from $5$, $N^1_{\downarrow}(5)$, all has out-degree 0 so $N^1_{\downarrow}(5)\subset\CI$. $\loverharpo\Gamma(K)$ inherits those directions from $\overharpo\Gamma(K_2)$ so $N(5)\setminus\{2,3\}\subset \CI$. Thus, we have proved 2(d). \\
        
        Assume condition (2), we want to show condition (1) by constructing such a pseudo-solvable group $G$. We assign 5 to $c$. Take a $3$-coloring $\{\CO, \CD, \CI\}$ of $\overline\Gamma$ where $\CD$ is the color that contains $5$, $\CI$ is the color that contains $N(5)\setminus\{a,b\}$, and $\CO$ is the remaining color. We direct the edges of $\overline\Gamma$ by assigning the arrows pointing from $\CO$ to $\CD$ to $\CI$. From the vertices $\{a,b\}$, we assign $2$ to the vertex that lies in $\CI$ and $3$ to the vertex that lies in $\CO$. Thus, $\overline\Gamma$ becomes a directed graph $\overharpo\Gamma$ with the assigned orientation. Denote $\overharpo F := \overharpo\Gamma\setminus\{\alpha,\beta,\gamma\}$, $n_o := |\CO\cap\overharpo F|, n_d:= |\CD\cap\overharpo F|, n_i := |\CI\cap \overharpo F|$. For any $v_j\in\CO\cap \overharpo F$, we equate $v_j$ with a prime $p_j$ such that the set $\CP = \{p_j\mid j = 1,\ldots, n_o\}$ is a set of distinct primes. Set $p$ to be the product of the distinct primes in $\CP$. By Dirichlet's theorem on arithmetic progressions, we pick a set of distinct primes $\CQ := \{q_k\mid k = 1,\ldots, n_d\}$ such that $q_k\equiv 1\pmod p$, for any $k$. We equate each vertex in $\CD\cap\overharpo F$ with a prime in $\CQ$. We define groups $P = C_{p_1}\times\cdots\times C_{p_{n_o}}$ and $Q = C_{q_1}\times\cdots\times C_{q_{n_d}}$. If $p_j\to q_k$ is an edge in $\overharpo\Gamma$, we let $C_{p_j}$ acts fixed point freely on $C_{q_k}$. If $p_j$ and $q_k$ are not adjacent to each other, we let $C_{p_j}$ acts trivially on $C_{q_k}$. Thus, we have defined a group action $P\acts Q$ and we denote the induced semidirect product $K = Q\rtimes P$. Since $K$ is solvable, it follows that a Hall $\{p_j, q_k\}$-subgroup of $K$ is a nontrivial product if and only if $p_j\text- q_k$ is an edge in $\overline\Gamma$. Thus, we know $\overharpo\Gamma(A_5\times K) \simeq \overharpo\Gamma[\CO\cup\CD\cup \{2\}]$. \\
        
        For any $v_j\in\CI$, let $N^1(v_j)$ be the set of primes in $\overharpo F$ with distance 1 to $v_j$ and $N^2(v_j)$ be the set of primes in $\overharpo F$ with distance 2 to $v_j$. If $N^1(v_j)$ is nonempty, let $B_{v_j}$ be a Hall $N^1(v_j)\cup N^2(v_j)$-subgroup of $K$. By the definition of $K$, we know $\Fit(B_{v_j})$ is a Hall $N^1(v_j)$-subgroup of $K$. Since $K$ has square-free order and so $B_{v_j}$ has square-free order. If $5\to v_j$ is an edge in $\overharpo\Gamma$, we consider the representation $t_4$ in Table~\ref{tab:rep of a5}. By arithmetic progression, we pick a prime $r_j$ such that $|A_5\times B_{v_j}|\mid r_j-1$. Let $r_j$ be large enough such that there exists an elementary abelian $r_j$-group $R_{v_j}$. $A_5$ acts on $R_{v_j}$ through the representation $t_4$ as described by Lemma~\ref{modular representation lemma}. By Lemma~\ref{representation lemma}, define $R_{v_j}$ to be a $A_5\times B_{v_j}$-representation. Therefore, if $f$ has an order that is divisible by $5$ or any primes in $N^1(v_j)$, then $f$ acts fixed point freely on $R_{v_j}$. If $5\to v_j$ is not an edge, we let $R_{v_j}$ be a $A_5\times B_{v_j}$-representation induced by the trivial representation $t_1$ of $A_5$. We let $\CR = \{r_j\mid j = 1,\ldots, n_i\}$. Based on the way the $r_j$ are chosen, by arithmetic progression, we can further require $\CR$ to be disjoint from $\CP\cup \CQ$. We let $J = R_{v_1}\times\cdots\times R_{v_{n_i}}$. We define $G = J\rtimes (K\times A_5)$ where the action $K\times A_5\acts J$ is induced by each $A_5\times B_{v_j}$ acts on $R_{v_j}$ as a group representation. By the way the actions are defined, we know the Frobenius digraph of the subgroup $R_{v_j}\rtimes (A_5\times B_{v_j})$ is isomorphic to the subgraph $\overharpo\Gamma[\{2,3,5,v_j\}\cup N^1(v_j)\cup N^2(v_j)]$. Therefore, $G$ is a pseudo-solvable group such that $\overline\Gamma(G)\simeq \overline\Gamma$. 
        \end{proof}

    \section{Category~(2): $\{p, 3, 5\}$ forms a triangle}
    
    We now discuss the second case. Assume that in the complement prime graph $\overline\Gamma(G)$, there exists a prime $p\notin\{2,3,5\}$ such that $\{p,3,5\}$ forms a triangle. We first show that all complement prime graphs of pseudo-solvable groups are 3-colorable. 
    \begin{theoremN}\label{3 colorable}
        Given a pseudo-solvable group $G$, $\overline\Gamma(G)$ is 3-colorable.
    \end{theoremN}
    \begin{proof}
        Given a pseudo-solvable group $G$, we look at the subgroup $K = N.(A_5\times H)$ where $\pi(G) = \pi(K)$. Let $K_1 = N.(D_{10}\times H)$ and $K_2 = N.(A_3\times H)$. If $3\mid |N|$ or $5\mid |N|$, then $\pi(K) = \pi(K_1)$ or $\pi(K) = \pi(K_2)$. Since $K_1$ and $K_2$ are solvable, $\overline\Gamma(K)$ is 3-colorable and so is $\overline\Gamma(G)$.\\
        
        Assume $(|N|, 15) = 1$. First note that if $\{2,3,5\}$ forms a triangle
        in $\overline\Gamma(G)$, then it is 3-colorable by Theorem \ref{235 triangle main theorem}.  Also, if $3-5$ is not an edge in $\overline\Gamma(G)$, then
        we are done by \cite[Theorem 6.6]{dawei2021prime}. So we may assume that $\{2,3,5\}$ is not a triangle, the 3 and 5 are linked by an edge in $\overline\Gamma(G)$.\\
        We give $\overline\Gamma(K)\setminus\{3\text -5\}$ the Frobenius digraph orientation based on $K_1, K_2$ discussed in Section 2, which we denoted $\overharpo\Gamma(K)$. Since we have at most one of the edges $3\to 2, 2\to 5$, by orienting the edge $3-5$ by $3\to 5$ we see that in $\{2,3,5\}$
        there is no directed 2-path.
        The directed graph deleting $\{2,3,5\}$ does not contain any directed 3-path. Therefore to show that the whole directed graph does not contain a directed 3-path, it suffices to show that all the edges starting from $2,3,5$ (which are by definition outgoing edges), end on vertices in $\CI$. This will imply that $\overharpo\Gamma(K)$ with the added  $3\to 5$ does not contain a directed $3$-path, which leads to the conclusion that $\overline\Gamma(G)$ is 3-colorable. Therefore, more formally we claim for any prime $p$ such that $q\to p$ is an edge for some $q\in \{2,3,5\}$, $p$ is a terminal vertex, i.e $p\in\CI$. Let $r\in \pi(G)\setminus\{p,2,3,5\}$ and $p\to r$ be an edge, we take a Hall $\{q,p,r\}$-subgroup $H_{qpr}$ and it must be a 2-Frobenius group of type $(r,p,q)$. Hence a Sylow $p$-subgroup $P$ is cyclic. Let $\pi  = \{p,2,3,5\}$. Consider a Hall $\pi$-subgroup of $K$. We take its chief series with $A_5$ as a top factor. Then, we take the last $p$-factor in the series and denote it $P_1$. We quotient everything before $P_1$ from $H_\pi$ and denote this new group $H_1$. $P_1$ is a chief factor so is an elementary abelian $p$-group. It is also a factor of the cyclic group $P$ so it is cyclic and thus $P_1\simeq C_p$. By the Schur-Zassenhaus theorem, $H_1 = P_1\rtimes (S.A_5)$ where $S$ is a 2-group. The automorphism group pf $P_1$ is cyclic, which implies that $A_5$ is a section of the kernel in the action of $S.A_5$ on $P_1$. This contradicts $q\to p$ being an edge and completes the proof of the theorem.
        
    \end{proof}
    The following lemma shows that category (1) is disjoint from category (2). 
    
    \begin{lemmaN}\label{2-3 2-5 not edges}
        Given $G$ a finite pseudo-solvable group, if the vertices $\{p,3,5\}$ form a triangle in $\overline\Gamma(G)$, where $p$ is some prime that is not 2, then $2\text-5$ and $2\text-3$ are not edges in $\overline\Gamma(G)$.
    \end{lemmaN}
    The proof will be similar to the proof of Lemma~\ref{3 not adjacent}, where we first prove some special cases and then generalize them.
    \begin{proof}
        Consider $K\simeq N.(A_5\times H)\leq G$ where $\pi(K) = \pi(G)$. If $3\mid |N|$ or $5\mid |N|$, $\overline\Gamma(K)$ cannot contain a triangle $\{p,3,5\}$ (which is seen in the usual way of going to suitable solvable subgroups). Thus $(|N|, 15) = 1$, and so by Lemma~\ref{3 not adjacent}, $3\text-p$ and $2\text-3$ cannot be edges a the same time. This proves that $2\text-3$ is not an edge in $\overline\Gamma(G)$.\\
        
        Now assume $2\text - 5$ is an edge in $\overline\Gamma(G)$. We will first show a special case $G = (P\rtimes Q).A_5$, where $P$ is an elementary abelian $p$-group and $Q$ is a 2-group. We assume $G$ has a normal series $$
        1\triangleleft P\triangleleft P\rtimes Q \triangleleft (P\rtimes Q).A_5.
        $$ Denote $N = P\rtimes Q$. By Schur-Zassenhaus theorem, since $P\triangleleft G$, $G$ has a subgroup $H\simeq Q.A_5$ such that $G= P\rtimes H$. Let $K_0$ be the kernel of the action $H\acts P$. Since $K_0$ centralizes $P$ and is normal in $H$, $K_0$ is a normal in $G= P\rtimes H$. Since $3\text-p, 5\text-p$ are edges in $\overline\Gamma(G)$, $K_0$ has no elements of order 3 or 5 so $K_0$ is a 2-group. Thus clearly $K_0\leq N$.
        
        WLOG, by working with $G/K_0\simeq P\rtimes (Q/K_0).A_5$, we assume $K_0=1$ and the action $H\acts P$ is faithful. Let $f\in G$ be a 5-element. Since $5\text-p, 5\text-2$ are edges in $\overline\Gamma(G)$, $\langle f\rangle$ acts on $N=P\rtimes Q$ frobeniusly. Since $N$ is a Frobenius kernel, $N$ is nilpotent so $N=P\times Q$. However, since the action of $H$ on $P$ is faithful, $Q = 1$. By Lemma~\ref{plms lemma}, $G=P.A_5$ so a 3-element $t$ acts on $P$ with a fixed point, implying $3-p$ is not an edge in $\overline\Gamma(G)$, which is a contradiction.\\
        
        As in the proof of Lemma~\ref{3 not adjacent}, we extend the special case to the general case. $G$ must be strictly pseudo-solvable so we take the subgroup $K\simeq N.(A_5\times H)$ where $(|N|, 15) = 1$. By working with a Hall subgroup, we assume $\pi(K) = \{2,3,5,p\}$. Take a chief series of $K$ where $A_5$ is the last factor. By quotienting out the subgroups before the last $p$-group factor, we assume $K$ has a normal series $1\triangleleft P\triangleleft P\rtimes Q\triangleleft (P\rtimes Q).A_5 = K$, where $P$ is an elementary abelian $p$-group, $Q$ is a 2-group. By the above special case, $2\text-5$ is not an edge in $\overline\Gamma(K)$ and thus the same statement holds in $\overline\Gamma(G)$. 
    \end{proof}
    Theorem~\ref{3 colorable} and Lemma~\ref{2-3 2-5 not edges} shows that category (1) and (2) describes two disjoint classes of prime graphs of pseudo-solvable groups. Together, they describe all the nontrivial prime graphs. 
    
    \begin{theoremN}\label{exact characterization}
        Let $G$ be a pseudo-solvable group and $\Gamma = \Gamma(G)$. Then, $\overline\Gamma$ satisfies exactly one of the following statements: \begin{enumerate}
            \item $\overline\Gamma$ is 3-colorable and triangle-free;
            \item $\overline\Gamma$ is 3-colorable and the $\{2,3,5\}$-subgraph is a triangle;
            \item $\overline\Gamma$ is 3-colorable and there exists a prime $p\in \pi(G)\setminus\{2,3,5\}$ such that the $\{p, 3,5\}$-subgraph is a triangle.
        \end{enumerate}
    \end{theoremN}
    \begin{proof}
        Combine Theorem~\ref{3 colorable} and Lemma~\ref{2-3 2-5 not edges}.
    \end{proof}
    
    It remains to classify all the $\overline\Gamma(G)$ where $\{p,3,5\}$ forms a triangle for some $p\neq 2$. By Lemma~\ref{dawei theorem 1}, all possible triangles of the complement prime graph $\overline\Gamma(G)$ have a base edge $3\text-5$. That condition restricts on how many copies of $A_5$ can appear in the chief factor of $G$. 
    
    \begin{lemmaN}\label{unique A5}
        Let G be a strictly pseudo-solvable group where $3\text-5$ is an edge in $\overline\Gamma(G)$, then $G$ has exactly one non-solvable chief factor which is isomorphic to $A_5$. 
    \end{lemmaN}
    \begin{proof}
        Let $\pi = \{2,3,5\}$. By definition, the nonsolvable sections of $G$ are isomorphic to $A_5$ so $G$ is $\pi$-separable. Let $H$ be a Hall $\pi$-subgroup of $G$. Let $N\triangleleft H$ be a solvable group such that all minimal normal subgroups of $H/N$ are nonsolvable. Denote $H_0 = H/N$ and let $M$ be a minimal normal subgroup of $H_0$. By the fact that a minimal normal subgroup of finite groups must be a direct product of copies of the same simple group, we know $M\simeq A_5^k$ for some $k\in \BZ^+$. If $k>1$, $M$ contains an element of order 15 and thus so does $G$, contradicting that $3\text-5$ is an edge. Thus, $M\simeq A_5$. $H_0$ acts on $M$ by conjugation so the action sends $H_0/C_{H_0}(M)$ to $\Aut(M)\simeq S_5$ injectively. Since $H_0$ does not contain an element of order $15$, $C_{H_0}(M)$ does not contain an element of order 3 or 5. Since $|H_0/C_{H_0}(M)|$ divides $|S_5| = 120$, $|H_0/C_{H_0}(M)|$ is divisible by $3$ and $5$ to the first power. Thus, $|H_0|$ is divisible by $3$ and $5$ only to the first power but not higher powers. By the group order, $M$ is the only section that is isomorphic to $A_5$ in $H_0$. Since $H_0 = H/N$ where $N$ is solvable, we know $H$ contains only one section of $A_5$ and so $G$ only contains one section of $A_5$. 
    \end{proof}
    Furthermore, we discover more properties of $\overline\Gamma(G)$. When $A_5$ acts on a $r$-group $R$, the action of $3$-elements is closely related to the actions of $5$-elements.
    \begin{lemmaN}\label{kmain}
	Let $G$ be a finite group and $N$ a normal subgroup of $G$ such that $G/N\cong A_5$ and $N$ is solvable.
    Suppose also that $|N|$ is divisible by a prime $p\not\in\{2,3,5\}$ such that $G$ does not contain any elements of orders $3p$ and $5p$.
    Now let $V$ be a finite $G$-module over the field with $q$ elements, where $q$ is a prime not dividing $|G|$. If the semidirect product $GV$ does not contain any elements of order $3q$, then $GV$ does not contain any elements of order $5q$.
    \end{lemmaN}
    The proof of Lemma~\ref{kmain} is long, so we present it in the Appendix. It follows immediately that the edge $3\text-r$ implies $5\text-r$ in the action of $A_5$ on $N$.
    
    \begin{lemmaN}\label{3p implies 5p}
    Let $G$ be a pseudo-solvable group. Suppose that $\overline\Gamma(G)$ contains a triangle with vertices $\{p, 3,5\}\subset\pi(G)$. Let $r\in\pi(G)$ be another prime, then $3\text-r$ implies that $5\text-r$ is also an edge in $\overline\Gamma(G)$.
    \end{lemmaN}
    \begin{proof}
    Let $H$ be a Hall $\{p,r,2,3,5\}$-subgroup of $G$. By Lemma~\ref{unique A5},  $H = N.A_5.M$ where $N, M$ are solvable. If $p\mid |M|$, then $3\text-p$ is an edge in $\overline\Gamma(A_5.M)$. Since $\Out(A_5)\simeq C_2$, $p$-elements in $M$ acts on $A_5$ trivially by conjugation, which produces elements of order $3p$, contradiction. Similarly, we know $r,3,5$ do not divide $|M|$ so $|M|$ is a $2$-group. We take a chief series of $N.A_5$ where $A_5$ is the top factor, $$
    1 = N_0\triangleleft N_1\triangleleft\cdots\triangleleft N_m = N\triangleleft N.A_5.$$ 
    For any $i$ such that the factor $N_{i+1}/N_i$ is an elementary abelian $r$-group, consider the group $M_i = (N.A_5)/N_i$ so it has an elementary abelian normal subgroup $Q_i = N_{i+1}/N_i$. Consider a Hall $\{2,3,5,p\}$-subgroup $P_i$ of $M_i$, so $P_i$ contains no elements of orders $3p$ or $5p$. Now consider the group action $P_i\acts Q_i$. By Lemma~\ref{kmain}, the group $Q_i\rtimes P_i$ has no elements of order $3r$ or $5r$. Therefore
    we see that $N.A_5$ has no elements of order $5r$. Since $M$ is a $2$-group, the same statement holds for $N.A_5.M$. Thus, $\overline\Gamma(G)$ has an edge $5\text-r.$ 
\end{proof}
    When $2\text-q$ is an edge in the complement prime graph, we next show that the Hall $\{2,3,5\}$-subgroup of $G$ has a strong connection with $\SL(2,5)$.
    
    \begin{lemmaN}\label{generalized quaternion lemma}
        Let $N$ be a $2$-group. If a Sylow $2$-subgroup of $N.A_5$ is generalized quaternion, then $N = C_2$ and $N.A_5\simeq \SL(2,5)$. 
    \end{lemmaN}
    \begin{proof}
        Identify a Sylow $2$-subgroup of $N.A_5$ with $Q_{2^n}$. If $n =3$, by calculating the order, we know $N = C_2$. The extension of $A_5$ by $C_2$ is either $C_2\times A_5$ or $\SL(2, 5)$, but a Sylow $2$-subgroup of $C_2\times A_5$ is not isomorphic to $Q_8$ so $N.A_5\simeq \SL(2,5)$ as desired.\\
        
        Assume $n>3$. The quotient $Q_{2^n}/N$ is isomorphic to a Sylow 2-subgroup of $A_5$, which is $V_4$. Since the abelianization $(Q_{2^n})^{ab}\simeq V_4$, we know $N$ is isomorphic to the commutator subgroup of $Q_{2^n}$, so $N$ is cyclic. Take a chief series of $N.A_5$, $
        1 = N_0\triangleleft N_1\triangleleft\cdots\triangleleft N_m = N\triangleleft N.A_5.$ Since a minimal normal subgroup of a finite group must be a direct product of the same simple group, the chief factors $N_{i+1}/N_i$ are elementary abelian 2-groups for $i\leq m-1$. The top factor is $N.A_5/N\simeq A_5$. Since $N$ is cyclic, its chief factors are also cyclic so $N_{i+1}/N_i\simeq C_2$ for $i\leq m-1$. Since $|N| = 2^m$, $N.V_4\simeq Q_{2^n}$ implies $m = n-2$. Since $N_1$ is normal in $N.A_5$ and also normal in $Q_{2^n}$, we know $N_1 = Z(Q_{2^n})$. Denote $M_i = (N.A_5)/N_i$ for $i\leq n-3$. For $i = 1$, we notice $M_1 = (N.A_5)/N_1$ so a Sylow $2$-subgroup of $M_1$ is $Q_{2^n}/Z(Q^{2^n}) = D_{2^{n-1}}$, the dihedral group of order $2^{n-1}$. Notice that $Z(D_{2^j})$ has order 2. Thus, for any index $1<j\leq n-3$, $N_i/N_{i-1} = Z(D_{2^{n-(i-1)}})$ so the Sylow 2-subgroup of $M_i$ is $D_{2^{n-(i-1)}}/Z(D_{2^{n-(i-1)}}) = D_{2^{n-i}}$. Therefore $M_{n-3} = (N.A_5)/N_{n-3}$ has Sylow 2-subgroups that are isomorphic to $D_8.$ However, $M_{n-3}$ is an extension of $A_5$ by $N_{n-2}/N_{n-3} = C_2$ which is either $C_2\times A_5$ or $\SL(2,5).$ Since neither of them has Sylow 2-subgroup isomorphic to $D_8$, this situation is impossible. 
    \end{proof}
    
    The following lemma tells us when $G$ has a subgroup that is isomorphic to $\SL(2,5)$. 
    \begin{lemmaN}\label{exists a sl25}
        Let $G = N.A_5$ where $N$ is a solvable group such that $(|N|, 15) = 1$. Assume $2\text-p$ is an edge in $\overline\Gamma(G)$, then a Hall $\{2,3,5\}$-subgroup of $G$ is isomorphic to $\SL(2,5)$.
    \end{lemmaN}
    \begin{proof}
        Since $G$ is \{2,3,5\}-separable, it has a Hall $\{2,3,5\}$-subgroup. Since $N$ is solvable, $G$ has a Hall $\{p,2,3,5\}$-subgroup $H$. Take a chief series of $G$, denoted $\Delta$, where $A_5$ is at the top: $0= N_0\triangleleft N_1\triangleleft\cdots\triangleleft N_m = N\triangleleft N.A_5$. Intersect the series with $H$ so we obtain a normal series $\Delta\cap H$ where each factor is a Hall $\{2,3,5,p\}$-subgroup of the original factor of $\Delta. $ Since $(|N|, 15) = 1$, all factors of $H\cap \Delta$ before $A_5$ are either elementary abelian 2-groups or elementary abelian $p$-groups. Now choose $i$ maximal such that $N_i/N_{i+1}$ is a $p-group$. The Sylow 2-subgroup of $G/N_i$ then acts frobeniusly on $N_i/N_{i+1}$ and thus is generalized quaternion. By Lemma~\ref{generalized quaternion lemma}, $G/N_i$ is isomorphic to $\SL(2,5)$, and clearly $N_i/N_{i+1}$ is abelian, but not cyclic. But if $j<i$ is chosen maximal such that is $N_i/N_{j}$ is a $p$-group and $N_i/N_{j+1}$ is not, then $N_i/N_{j}$ cannot act frobeniusly on $N_j/N_{j+1}$ which contradicts $2\text-p$ being an edge in $\overline\Gamma(G)$. Hence 2 does not divide $|N_{i+1}|$ and the proof is complete.

\end{proof}

\begin{figure}[ht!]
    \centering
        \subfigure[$\rho_1$]{
            \begin{tikzpicture}[scale = 1]
        \node (2) at (-1,-1) {$2$};
        \node (3) at (-1,1) {$3$};
        \node (5) at (1,1) {$5$};
        \node (p) at (1,-1) {$r$};
        \path [-] (3) edge node {} (5);
            \end{tikzpicture} 
            }\;\;\;\;\;\;
           \subfigure[$\rho_2$]{ \begin{tikzpicture}[scale = 1]
        \node (2) at (-1,-1) {$2$};
        \node (3) at (-1,1) {$3$};
        \node (5) at (1,1) {$5$};
        \node (p) at (1,-1) {$r$};
        \path [-] (3) edge node {} (5);
        \path [-] (3) edge node {} (p);
        \path [-] (5) edge node {} (p);
        \path [-] (2) edge node {} (p);
            \end{tikzpicture}
            }\;\;\;\;\;\;
            \subfigure[$\rho_6$]{ \begin{tikzpicture}[scale = 1]
        \node (2) at (-1,-1) {$2$};
        \node (3) at (-1,1) {$3$};
        \node (5) at (1,1) {$5$};
        \node (p) at (1,-1) {$r$};
        \path [-] (3) edge node {} (5);
        \path [-] (5) edge node {} (p);
            \end{tikzpicture}
            }\\
            \subfigure[$\rho_7$]{ \begin{tikzpicture}[scale = 1]
        \node (2) at (-1,-1) {$2$};
        \node (3) at (-1,1) {$3$};
        \node (5) at (1,1) {$5$};
        \node (p) at (1,-1) {$r$};
        \path [-] (3) edge node {} (5);
        \path [-] (5) edge node {} (p);
        \path [-] (2) edge node {} (p);
            \end{tikzpicture}
            }\;\;\;\;\;\;
            \subfigure[$\rho_9$]{ \begin{tikzpicture}[scale = 1]
        \node (2) at (-1,-1) {$2$};
        \node (3) at (-1,1) {$3$};
        \node (5) at (1,1) {$5$};
        \node (p) at (1,-1) {$r$};
        \path [-] (3) edge node {} (5);
        \path [-] (2) edge node {} (p);
            \end{tikzpicture}
            }
    \caption{Subgraphs constructed by the representations of $\SL(2,5)$}
    \label{fig: possible graphs of case 2}
    \end{figure}
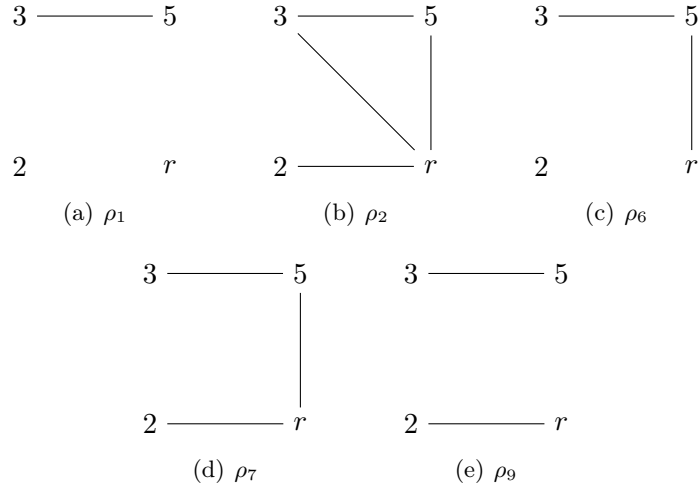

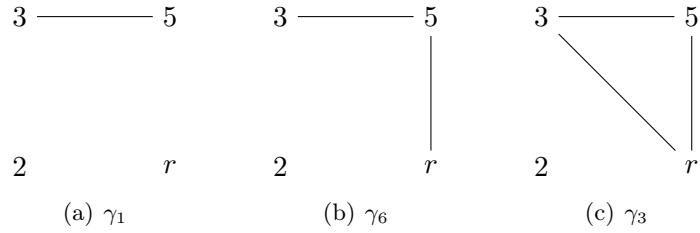
\begin{figure}[ht!]
    \centering\subfigure[$\gamma_1$]{
            \begin{tikzpicture}[scale = 1]
        \node (2) at (-1,-1) {$2$};
        \node (3) at (-1,1) {$3$};
        \node (5) at (1,1) {$5$};
        \node (p) at (1,-1) {$r$};
        \path [-] (3) edge node {} (5);
            \end{tikzpicture} 
            }\;\;\;\;\;\;
            \subfigure[$\gamma_6$]{ \begin{tikzpicture}[scale = 1]
        \node (2) at (-1,-1) {$2$};
        \node (3) at (-1,1) {$3$};
        \node (5) at (1,1) {$5$};
        \node (p) at (1,-1) {$r$};
        \path [-] (3) edge node {} (5);
        \path [-] (5) edge node {} (p);
            \end{tikzpicture}
            }\;\;\;\;\;\;
          \subfigure[$\gamma_3$]{ \begin{tikzpicture}[scale = 1]
        \node (2) at (-1,-1) {$2$};
        \node (3) at (-1,1) {$3$};
        \node (5) at (1,1) {$5$};
        \node (p) at (1,-1) {$r$};
        \path [-] (3) edge node {} (5);
        \path [-] (3) edge node {} (p);
        \path [-] (5) edge node {} (p);
            \end{tikzpicture} 
            }
    \caption{Subgraphs constructed by the representations of $C_2.S_5$}
    \label{fig: conditioned graphs of case 2}
    \end{figure}

\begin{figure}[ht!]
    \centering
        \subfigure[edges $3\text-5$, $3\text-r$]{
            \begin{tikzpicture}[scale = 1]
        \node (2) at (-1,-1) {$2$};
        \node (3) at (-1,1) {$3$};
        \node (5) at (1,1) {$5$};
        \node (p) at (1,-1) {$r$};
        \path [-] (3) edge node {} (5);
        \path [-] (3) edge node {} (p);
            \end{tikzpicture} 
            }\;\;\;\;\;\;
            \subfigure[edges $3\text-5$, $3\text-r$, $2\text-r$]{ \begin{tikzpicture}[scale = 1]
        \node (2) at (-1,-1) {$2$};
        \node (3) at (-1,1) {$3$};
        \node (5) at (1,1) {$5$};
        \node (p) at (1,-1) {$r$};
        \path [-] (3) edge node {} (5);
        \path [-] (3) edge node {} (p);
        \path [-] (2) edge node {} (p);
            \end{tikzpicture} 
            }
    \caption{Impossible subgraphs}
    \label{fig: unknown graphs of case 2}
    \end{figure}
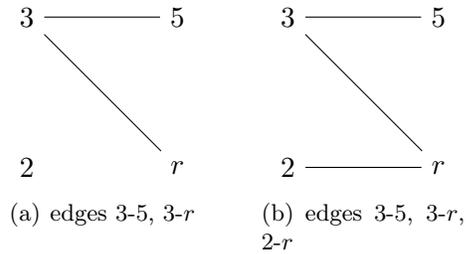

    Suppose that $G$ is a group such that $\overline\Gamma(G)$ where $\{p,3,5\}$ forms a triangle for some $p\neq 2$. (This is the situation that we are interested in right now.) For an arbitrary prime $r\in \pi(G)\setminus\{2,3,5\}$, Figures~\ref{fig: possible graphs of case 2} and ~\ref{fig: conditioned graphs of case 2} list all the possible edge relations between the vertices $\{r,2,3,5\}$,
    as will be shown below. Lemma~\ref{3p implies 5p} shows that the subgraphs in Figure~\ref{fig: unknown graphs of case 2} do not occur. Subgraphs in Figure~\ref{fig: conditioned graphs of case 2} are different from subgraphs in Figure~\ref{fig: possible graphs of case 2} by the edge $2\text-r$. Lemma~\ref{figure 1 or figure 3} will establish this important difference where Lemma~\ref{exists a sl25} and $\SL(2,5)$ will play a crucial role. 
    (For the labels of the representations, see the tables below.)
    
    \begin{lemmaN}\label{bringing C2 up}
        Let $G$ be a pseudo-solvable group where $\{p,3,5\}$ forms a triangle in $\overline\Gamma(G)$. Assume there exists a prime $q$ such that 2-q is an edge in $\overline\Gamma(G)$. Let $r$ be a prime divisor of $G$ that is not 2,3, or 5. Let $H$ be a Hall $\{r,2,3,5\}$-subgroup of $G$. Then, $H \simeq R.Q.A_5.M$ where $R$ is a Sylow $r$-subgroup of $H$, $Q$ is a cyclic of order 2, and $M$ is a 2-group. 
    \end{lemmaN}
    \begin{proof}
        Let $K$ be a Hall $\{p,r,2,3,5\}$-subgroup of $G$. By Lemma~\ref{unique A5}, $K = N_0.A_5.M_0$ where $N_0, M_0$ are solvable. Since $\{p,3,5\}$ induces a triangle in $\overline\Gamma(K)$, it is easy to see that 3 and 5 do not divide $|N_0|$. By Lemma~\ref{exists a sl25}, a Hall $\{2,3,5\}$-subgroup of $N_0.A_5$ is isomorphic to $\SL(2,5)$. Let $H$ be a Hall $\{r,2,3,5\}$-subgroup of $G$ and $H \simeq N.A_5.M$ where $N, M$ are solvable and $(|N|, 15)=1$. As a consequence, the Hall $\{2,3,5\}$-subgroups of $N.A_5$ are isomorphic to $\SL(2,5)$ which implies that the Sylow 2-subgroups of $N$ are cyclic of order 2. Assume $H$ has a normal series $1 < P_1 \leq P_1.C_2\leq P_1.C_2.P_2 = N\leq N.A_5\leq N.A_5.M = H$ for some $r$-subgroups $P_1$, $P_2$. By the Schur-Zassenhaus theorem, the quotient $C_2.P_2$ is a semidirect product and hence even a direct product. Thus, after denoting $R:= P_1.P_2$, we may modify the normal series as $1 < R\leq R.C_2\leq R.C_2.A_5\leq R.C_2.A_5.M$. Since $3\text-5$ is an edge in $\overline\Gamma(A_5.M)$ and $\Out(A_5)\simeq C_2$, we know the order of $M$ is coprime to 15, and any $r$-subgroup of $M$ has to act trivially on $C_2.A_5$
        and thus can be "moved down" into $R$, and so $M$ is a 2-group and the assertion follows.
    \end{proof}

Recall from Lemma \ref{plms lemma} that if $G$ is a finite group with normal subgroup $H$ and $G/H\simeq A_5$, and $C_H(t)=1$ for some $3$-element $t$, then $H$ is a $2$-group. This
will be useful in the proof of the following lemma. 
    
    \begin{lemmaN}\label{figure 1 or figure 3}
        Let $G$ be a pseudo-solvable group where $\{p,3,5\}$ induces a triangle in $\overline\Gamma(G)$. Let $r$ be a prime divisor of $\pi(G)$. Then $\{r,2,3,5\}$-subgraphs are not graphs in Figure~\ref{fig: unknown graphs of case 2}. Moreover, if there is at least one $\{r,2,3,5\}$-subgraph in $\overline\Gamma(G)$ as in Figure~\ref{fig: possible graphs of case 2} (b), (d), or (e), then there is no $\{r,2,3,5\}$-subgraph in $\overline\Gamma(G)$ as in 
        Figure~\ref{fig: conditioned graphs of case 2} (c). 
        
    \end{lemmaN}
    \begin{proof}
        The first statement, that subgraphs in Figure~\ref{fig: unknown graphs of case 2} do not occur, follows immediately from Lemma~\ref{3p implies 5p}.
    
        To prove the second statement, suppose that $\overline\Gamma(G)$ has a $\{q,2,3,5\}$-subgraph that is of a kind in Figure~\ref{fig: possible graphs of case 2} (b),(d), or (e), and assume that a $\{r,2,3,5\}$-subgraph as in Figure~\ref{fig: conditioned graphs of case 2}(c) also appears. (We need to obtain a contradiction.) In particular, $2\text-q$ is an edge. By Lemma~\ref{exists a sl25} and Lemma~\ref{bringing C2 up}, a Hall $\{r,2,3,5\}$-subgroup $H$ has a normal series $1 < R \leq N\leq M\leq H$ where $R$ is a Sylow $r$-subgroup, $M/R$ is isomorphic to $\SL(2,5)$ and $N/R\simeq C_2$, $M/N\simeq A_5$, and $H/M$ is a  2-group. Since $2\text-q$ is an edge, the Sylow 2-subgroups of $G$ are generalized quaternion. We let $T$ be a Sylow $2$-subgroup of $H$. We assume $R$ is elementary abelian up to quotienting $r$-groups from $H$, where Figure~\ref{fig: conditioned graphs of case 2}(c) is preserved. By Schur-Zassenhaus theorem, $R$ has a complement $K\leq H$ such that $R\rtimes K\simeq H$. Since $N/R$ is normal in $H/R$ and is of order 2, $N/R$ is contained in the center of $K$. Write $Z$ for this central subgroup of order 2 of $K$. By Maschke's theorem, $R$ is a completely reducible $K$-module so $R = \bigoplus_{i=1}^m V_i$ where $V_i$ are irreducible $K$-modules. If $Z$ acts Frobeniusly on each $V_i$, then a Sylow 2-subgroup $T$ of $H$ acts Frobeniusly on $R$. It implies $2\text- r$ is an edge in $\overline\Gamma(H)$, which contradicts that Figure~\ref{fig: conditioned graphs of case 2}(c) is a subgraph. Therefore, WLOG, assume $Z$ acting on $V_1$ with a fixed point. By Zassenhaus decomposition, $V_1\simeq [V_1, Z]\oplus C_{V_1}(Z)$ where $Z$ acts Frobeniusly on $[V_1, Z]$ and trivially on $C_{V_1}(Z)$. However, since $V_1$ is an irreducible $K$-module and $Z$ does not act on $V_1$ Frobeniusly, we know $V_1 = C_{V_1}(Z)$ which means the action is trivial. We rewrite a normal series of $H$, $1\leq X< Y_1 < Y_2 < M\leq H$ where $X\simeq V_2\oplus\cdots\oplus V_m$, $Y_1/X\simeq C_2$, $Y_2/Y_1\simeq V_1$, and $M/Y_2\simeq A_5$. Consider the group $M/Y_1\simeq V_1.A_5$. Since $3\text-r$ is an edge in Figure~\ref{fig: conditioned graphs of case 2}(c), any $3$-element $t$ acts Frobeniusly on $Y_2/Y_1\simeq V_1$. By Lemma~\ref{plms lemma}, $V_1$ must be a $2$-group, contradicting that $V_1$ is a $r$-group. Therefore, we conclude that Figure~\ref{fig: conditioned graphs of case 2}(c) is impossible to be the $\{r, 2,3 ,5\}$-subgraph of $\overline\Gamma(G)$.
    \end{proof}
    
    \begin{table}[h!]
    \centering
    \begin{tabular}{|| c | c | c c c c c c c c c||}
    \hline
    Index & Degrees & \multicolumn{9}{c||}{Orders}\\
    \hline
     & & 1 & 10 & 10 & 2 & 5 & 5 & 3 & 6 & 4\\
    \hline
    $\rho_1$ & 1 & 0 & 0 & 0 & 0 & 0 & 0 & 0 & 0 & 0\\
    \hline
    $\rho_2$ & 2 &0 & - & - & - & - & - & - & - & -\\ 
    \hline
    $\rho_3$ & 2 & 0 & - & - & - & - & - & - & - & -\\
    \hline
    $\rho_4$ & 3 & 0 & 0 & 0 & 0 & 0 & 0 & 0 & 0 & 0\\
    \hline
    $\rho_5$ & 3 & 0 & 0 & 0 & 0 & 0 & 0 & 0 & 0 & 0\\
    \hline
    $\rho_6$ & 4 & 0 & - & - & 0 & - & - & 0 & 0 & 0\\
    \hline
    $\rho_7$ & 4 & 0 & - & - & - & - & - & 0 & - & -\\
    \hline
    $\rho_8$ & 5 & 0 & 0 & 0 & 0 & 0 & 0 & 0 & 0 & 0\\
    \hline
    $\rho_9$ & 6 & 0 & - & - & - & 0 & 0 & 0 & - & -\\
    \hline
    \end{tabular}
    \caption{The values of $\det(\rho(g)-I)$ for irreducible representations and conjugacy class representatives of $\SL(2,5)$}  \label{tab: new rep of sl25}  
    \end{table}
    
    \begin{table}[h!]
        \centering
        \begin{tabular}{||c| c| c c c c c c c c c c c c||}
            \hline
            Index & Degrees & \multicolumn{12}{c||}{Orders}\\
            \hline
             & & 1 & 6 & 2 & 3 & 6 & 5 & 4 & 8 & 10 & 8 & 6 & 2\\
            \hline
            $\gamma_1$ & 1 & 0 & 0 & 0 & 0 & 0 & 0 & 0 & 0 & 0 & 0 & 0 & 0\\
            \hline
            $\gamma_2$ & 1 & 0 & - & - & 0 & - & 0 & 0 & - & 0 & - & 0 & 0\\
            \hline
            $\gamma_3$ & 4 &0 & - & 0 & - & - & - & - & - & - & - & - & -\\
            \hline
            $\gamma_4$ & 4 & 0 & 0 & 0 & 0 & 0 & - & 0 & 0 & - & 0 & 0 & 0\\
            \hline
            $\gamma_5$ & 4 & 0 & 0 & 0 & 0 & 0 & - & 0 & 0 & - & 0 & 0 & 0\\
            \hline
            $\gamma_6$ & 4 & 0 & 0 & 0 & 0 & 0 & - & - & - & - & - & - & -\\
            \hline
            $\gamma_7$ & 4 & 0 & 0 & 0 & 0 & 0 & - & - & - & - & - & - & -\\
            \hline
            $\gamma_8$ & 5 & 0 & - & 0 & 0 & - & 0 & 0 & 0 & 0 & 0 & 0 & 0\\
            \hline
            $\gamma_9$ & 5 & 0 & 0 & 0 & 0 & 0 & 0 & 0 & 0 & 0 & 0 & 0 & 0\\
            \hline
            $\gamma_{10}$ & 6 & 0 & 0 & 0 & 0 & 0 & 0 & 0 & 0 & 0 & 0 & 0 & 0\\
            \hline
            $\gamma_{11}$ & 6 & 0 & 0 & 0 & 0 & 0 & 0 & - & - & - & - & - & -\\
            \hline
            $\gamma_{12}$ & 6 & 0 & 0 & 0 & 0 & 0 & 0 & - & - & - & - & - & -\\
            \hline
        \end{tabular}
        \caption{The values of $\det(\rho(g)-I)$ for irreducible representations and conjugacy class representatives of $C_2.S_5$}
        \label{tab:rep of c2a5}
    \end{table}
    
    For the backward direction, we use the irreducible representations of $\SL(2,5)$ and $C_2.S_5$, with GAP index SmallGroup(120,90), to construct the subgraphs in Figure~\ref{fig: possible graphs of case 2} and Figure~\ref{fig: conditioned graphs of case 2}, respectively. As we constructed subgraphs in Figure~\ref{fig: possible graphs of case 1} using $A_5$, we examine the fixed points of the irreducible representations with GAP~\cite{gap}. Using the same convention as Table~\ref{tab:rep of a5}, Table~\ref{tab: new rep of sl25} and Table~\ref{tab:rep of c2a5} describe fixed points of the irreducible representations of $\SL(2,5)$ and $C_2.S_5$ respectively. By Lemma~\ref{modular representation lemma}, $\rho_1$ of $A_5$ gives Figure~\ref{fig: possible graphs of case 2}(a); $\rho_2$ gives Figure~\ref{fig: possible graphs of case 2}(b); $\rho_6$ gives Figure~\ref{fig: possible graphs of case 2}(c); $\rho_7$ gives Figure~\ref{fig: possible graphs of case 2}(d); $\rho_9$ gives Figure~\ref{fig: possible graphs of case 2}(e). For $C_2.A_5$, $\gamma_1$ gives Figure~\ref{fig: conditioned graphs of case 2}(a); $\gamma_6$ gives Figure~\ref{fig: conditioned graphs of case 2}(b); $\gamma_3$ gives Figure~\ref{fig: conditioned graphs of case 2}(c). Combining Lemma~\ref{figure 1 or figure 3}, we give the following theorem. 
    
    \begin{theoremN}\label{category 2 description}
        Let $\Gamma$ be a simple graph, then the following are equivalent:
        \begin{enumerate}
            \item $\Gamma\simeq \Gamma(G)$ for some pseudo-solvable group $G$ such that $\{p, 3,5\}$ induces a triangle in $\overline\Gamma(G)$ for some prime $p\ne 2$;
            \item $\overline\Gamma$ with a choice of vertices $\{\alpha, \beta, \gamma\}\subset V(\Gamma)$ satisfies the following conditions:
            \begin{enumerate}
                \item $\overline\Gamma$ is 3-colorable;
                \item $\alpha\text-\beta$ is an edge for all triangles in $\overline\Gamma$, and at least one such triangle exists;
                \item $\overline\Gamma$ has a choice of coloring $\CO, \CD, \CI$ such that $N(\{\alpha, \beta, \gamma\})\setminus\{\alpha,\beta,\gamma\}\subset \CI$;
                \item The set $\mathcal{S} := \{\{q, \alpha, \beta, \gamma\}\text{-subgraph}\subset\overline\Gamma: q\in V(\overline\Gamma)\setminus \{\alpha, \beta, \gamma\}\}$ contains only graphs  from Figure~\ref{fig: possible graphs of case 2}, or only graphs from Figure~\ref{fig: conditioned graphs of case 2}. 
            \end{enumerate}
        \end{enumerate}
    \end{theoremN}
    \begin{proof} We show the forward direction and then the backward. For the forward direction, assume Condition~(1). Let $\alpha=3$, $\beta=5$, and $\gamma=2$.
    Theorem~\ref{3 colorable} implies (2a), and the proof of Theorem~\ref{3 colorable} implies (2c). Lemma~\ref{dawei theorem 1} implies (2b). It remains to show (2d). First note that by Lemma \ref{2-3 2-5 not edges} we know that the subgraph of $\overline\Gamma(G)$ induced by $\{2,3,5\}$ is fully determined. So if $r$ is a prime with
    $r\in \pi(G)\backslash \{2,3,5\}$, let $\pi_r = \{r,2,3,5,\}$ and $H_r$ a Hall $\pi_r$-subgroup of $G$. 
    Write $\Gamma_r=\overline\Gamma(H)$. Since $\Gamma_r[\{2,3,5\}]$ is already fully determined, there
    are theoretically only eight possibilities for $\Gamma_r$,
    depending on whether there are edges between 2,3,5 and $r$ or not. Two of these, namely the ones in Figure~\ref{fig: unknown graphs of case 2}, cannot occur by Lemma~\ref{figure 1 or figure 3}, reducing the number of possibilities to 6.
    Those six graphs are listed in Figures ~\ref{fig: possible graphs of case 2} and 3, but two of them are listed twice (namely those which have no edge between 2 and $r$.\\
    To complete the proof of (2d), we must show that either $\Gamma_s$ is in Figure ~\ref{fig: possible graphs of case 2} for all $s\in \pi(G)\backslash \{2,3,5\}$, or $\Gamma_s$ is in Figure 3 for all $s\in \pi(G)\backslash \{2,3,5\}$. Now both Figure 2 and Figure 3 list all (two) possible subgraphs $\Gamma_r$ where there is no edge between 2 and $r$, so those do not pose a problem. Now if $\overline\Gamma(G)$ has at least one subgraph
    $\Gamma_r$ of Figure 2 such that $2-r$ is an edge in $\Gamma_r$ then by Lemma~\ref{figure 1 or figure 3}
    we conclude that all subgraphs $\Gamma_s$ ($s\in \pi(G)\backslash \{2,3,5\}$) are in Figure~\ref{fig: possible graphs of case 2}. Likewise, if $\overline\Gamma(G)$ has at least one subgraph
    $\Gamma_r$ of Figure~\ref{fig: conditioned graphs of case 2} which is the unique subgraph in that figure such that $2-r$ is an edge in $\Gamma_r$, then again by Lemma~\ref{figure 1 or figure 3} it follows that all subgraphs $\Gamma_s$ ($s\in \pi(G)\backslash \{2,3,5\}$) are in Figure 3. This completes the proof of (2d).

        We now show the backward direction. Let $\Gamma$ be a simple graph such that its complement
        $\overline\Gamma$ satisfies Condition (2). $\overline\Gamma$ can be decomposed into two parts, the $\{\alpha, \beta,\gamma\}$-subgraph and the subgraph $\Delta$ generated by the other vertices. $\Delta$ is 3-colorable and triangle-free. We give $\Delta$ a coloring $\CO, \CD, \CI$ which satisfies (2c). We identify the vertices $\gamma =2$ and $\{\alpha, \beta\}=\{3,5\}$ in such a way as to be consistent with (2d).
        We first orientate $\Delta$ by directing edges pointing from $\CO$ to $\CD$ then to $\CI$, denote this graph $\overharpo\Delta$. We repeat the construction that we did in Theorem~\ref{235 triangle main theorem} and use the notation introduced there. That is, we identify vertices in $\CO$ with distinct primes $\CP$ and then identify vertices in $\CD$ with distinct primes $\CQ$ such that $Q\rtimes P$ has prime graph isomorphic to $\Delta\cap(\CO\cup\CD)$. Let $S = \SL(2,5)$ or $C_2.S_5$ depending on $\mathcal S$ having all its elements in  Figure~\ref{fig: possible graphs of case 2} or Figure~\ref{fig: conditioned graphs of case 2}, respectively. (If all elements of $\mathcal S$ are in Figure 2 as well as in Figure 3, then the the choice does not matter.) $S$ takes the place of $A_5$ in the Theorem~\ref{235 triangle main theorem} construction. By Lemma~\ref{representation lemma}, for any vertex $v_j\in \CI$, we assign to $v_j$ a sufficiently large prime $r_j$ and a group $R_{v_j}$. Then we let $S\times (Q\rtimes P)$ act on $R_{v_j}$ in the following way. For any prime $s\in \CP \cup \CQ$ the group $C_s$ acts on $R_{v_j}$ Frobeniusly if and only if $s \to r$ is an edge in $\overharpo\Delta$. 
        $S$ acts on $R_{v_j}$ by the representation that the $\{v_j, \alpha, \beta, \gamma\}$-subgraph corresponds in Figure~\ref{fig: possible graphs of case 2} or Figure~\ref{fig: conditioned graphs of case 2}. Therefore, we have a group action $S\times (Q\rtimes P)$ on $R$ which has a prime graph that is isomorphic to $\Gamma$ through the vertex assignment we have described. The corresponding semidirect product $G = R\rtimes (S\times (Q\rtimes P) )$ is the desired group. 
    \end{proof}
 
 \textup{   
Finally, Theorems~\ref{235 triangle main theorem} and ~\ref{category 2 description} together provide the desired classification of pseudo-solvable groups by their prime graphs, which was the main goal of the paper. We mention an interesting and a little less technical consequence of our classification.
}\\
    
    \begin{corollaryN}\label{cor:complete classification}
        Let $G$ be a pseudo-solvable group and $\Gamma = \Gamma(G)$, then $\Gamma$ is in exactly one of the following categories.
        \begin{enumerate}
            \item $\overline\Gamma$ is 3-colorable and triangle-free;
            \item $\overline\Gamma$ is 3-colorable and the $\{2,3,5\}$-subgraph is a triangle. For an arbitrary prime divisor $r$, all the $\{r,2,3,5\}$-subgraphs are of types in Figure~\ref{fig: possible graphs of case 1};
            \item $\overline\Gamma$ is 3-colorable and there exists an odd prime $p$ such that the $\{p,3,5\}$-subgraph is a triangle. For an arbitrary prime divisor $r$, all the $\{r,2,3,5\}$-subgraphs are of types in Figure~\ref{fig: possible graphs of case 2};
            \item $\overline\Gamma$ is 3-colorable and there exists an odd prime $p$ such that the $\{p,3,5\}$-subgraph is a triangle. For an arbitrary prime divisor $r$, all the $\{r,2,3,5\}$-subgraphs are of types in Figure~\ref{fig: conditioned graphs of case 2}.
        \end{enumerate}
    \end{corollaryN}

    \section{Outlook}

    Non-abelian simple groups whose order is divisible by exactly three different primes are often called $K_3$-groups. There are eight $K_3$ groups, and $A_5$ is one of them. If we call pseudo-solvable groups "pseudo $A_5$-solvable" groups, then it is a natural problem to study the prime graphs of pseudo $T$-solvable groups for other nonabelian simple groups $T$. In fact, in \cite{edwardskellerpesaklatha} the authors completely classify
    the prime graphs of pseudo $T$-solvable groups for all the remaining $K_3$ groups. It turns out that for all of these only the case similar to the case that there is a $\{2,3,5\}$ triangle for $A_5$ can occur, whereas the second case for $A_5$, that there
    is a triangle $\{3,5,p\}$ for some odd prime $p$ and which allows for an arbitrary number of such triangles, is unique to $A_5$.\\
    For further study One can now, of course, look into pseudo $T$-solvable groups for yet other nonsolvable groups $T$, such as the Suzuki groups. It might also be interesting to consider other graph theory aspects of prime graphs of pseudo-solvable groups, such as connectivity or a yet to be developed notion of the minimality of the graphs. 

    
    \section{Acknowledgement}
    This research was conducted at Texas State University under NSF-REU grant DMS-1757233 and NSA grant H98230-21-1-0333 during the summer of 2021. The authors thank NSF and NSA for the financial support. The first, third, fourth, and fifth thank Texas State University for running the REU online during this difficult period of social distancing and providing a welcoming and supportive work environment. Those authors also thank their mentor, the second author Dr. Thomas Michael Keller, for his invaluable advice and guidance throughout this project. The sixth author, the director of the REU program, is recognized for conducting an inspired and successful research program. The results in this paper were mainly discovered by the first author under the guidance of the second author, while the other members worked out another paper. 
    \\

\section{Appendix}
We present a proof for Lemma~\ref{kmain} in this section. 

\begin{proof}
Suppose the result is not true and  let $G$, $V$ be a minimal counterexample, i.e., $|G|+|V|$ is minimal.\\

Our first goal is to show that $|N|$ is not divisible by 3 and by 5. Assume that 3 divides $|N|$. Then there exist normal subgroups $S$ and $T$ of $G$ with $T<S\leq N$ such that $S/T$ is a chief factor of $G$ of order divisible by 3, and we may also choose $S,T$ such that  3 does not divide $|N/S|$. Then $p$ does not divide $|N/T|$ because otherwise any Hall $\{3,p\}$-subgroup of $G/T$ (which clearly exists) would be a 2-Frobenius group, making the Sylow $p$-subgroup of $G/T$ cyclic, but since there are 3-elements and 5-elements in $G/T$ which act frobeniusly on chief factors of $p$-power order of $N/T$ (which now must be of order $p$) and those 3- and 5-elements generate a nonabelian group, we get a contradiction to the automorphism group of cyclic groups being abelian. So $p$ indeed does not divide $|N/T|$. Therefore $S/T$ is a section of a 3-subgroup of $G$ acting frobeniusly on a chief factor of $p$-power order of $G$ in $T$, and hence that 3-subgroup must be cyclic. This shows that $S/T$ is cyclic of order 3.  (We remark that exploiting the Frobenius action of the Sylow 3-subgroups of $G$ on $V$ would yield $|S/T|=3$ much faster, but we want to imply the same argument also for the prime 5 for which we do not have the fixed point free action of Sylow 5-subgroups on $V$, and hence we only used the nonexistence of elements of order $3p$ here, since there are also no elements of order $5p$.)\\

Now $(G/T)/C_{G/T}(S/T)$ has order 1 or 2. If the order is 2, then $G/C_G(S/T)$ is of order 2 and so  $C_G(S/T)$ acting on $V$ will be a smaller counterexample, contradicting our minimal choice of $|G|$, $|V|$. Thus the order is 1, meaning that $S/T$ is a central chief factor, and then $N/T$ is 3-nilpotent, and $G$ has a normal subgroup $R$ such that $G/R$ is a central extension of $A_5$ by $C_3$ (the cyclic group of order 3). Since the Schur multiplier of $A_5$ is 2, it is clear that this extension must be split, and hence $G$ has a normal subgroup $R_1$ such that $G/R_1$ is of order 3. Then $R_1$, $V$ is a smaller counterexample to our theorem, which is a contradiction. This shows that the assumption of 3 dividing $|N|$ must be refuted, so 3 does not divide $|N|$.\\

A similar argument shows that 5 does not divide the order of $N$: Here one will, otherwise, find a factor group of order 1 or 2 or 4 and can handle the cases 2 and 4 by induction. Also, once one has found a normal subgroup $T_1$ in $G$ of index 5 in $G$ then $T_1$, $V$ is no longer a counterexample, so there are no elements of order $5q$ in $T_1V$, and since $G/T_1$ is part of a cyclic group acting frobeniusly on some $p$-chief factor of $G$, then clearly  there cannot exist any elements of order $5q$ in $GV$, contradicting $G$, $V$ being a counterexample. Therefore 5 does not divide $|N|$.\\

We remark that since $|N|$ is not divisible by 3, then for any Sylow 3-subgroup $R$ of $G$ we have that $R$ acts coprimely on $N$, and so $N$ has an $R$-invariant Sylow $p$-subgroup $P$. Then $PV$ is a subgroup of $GV$ on which $R$ acts frobeniusly, and thus $PV$ is nilpotent. This shows that $P\leq C_G(V)$.\\

Since $(|G|,|V|)=1$, we know that $V$ is a completely reducible $G$-module, and since for any irreducible $G$-submodule $W$ of $V$ clearly $G/C_G(W)$ is divisible by 3 by our hypothesis that $GV$ does not contain any elements of order $3q$, we see that the action of $G$ on $W$ satisfies the hypothesis of the theorem, and if $W<V$, then $G$, $W$ is not a counterexample and thus all elements of order 5 act fixed point freely on $W$. Since $W$ was arbitrary, then it follows that all elements of order 5 act fixed point freely on $V$, contradicting $G$ being a counterexample. This shows that $V$ is irreducible as $G$-module.\\

Now put $\pi=\{2,3,5\}$ and suppose that $|G/C_G(V)|_{\pi}=60$. By the Schur-Zassenhaus theorem, $G/C_G(V)$ has a Hall $\pi$-subgroup isomorphic to $A_5$ which acts on $V$. Then from Lemma \ref{plms lemma} 
we know that the Sylow 3-subgroups of $G/C_G(V)$ will have nontrivial fixed points on $V$, contradicting our hypothesis that $GV$ does not contain any elements of order $3q$. Therefore $|G|_{\pi}>60$, and since 3 and 5 both do not divide $|N|$, it follows that 2 divides $|N|$. \\

Next we observe that if $1=N_1<N_2<\dots < N_k=G$ is a chief series of $G$, then necessarily $A_5$ is the top chief factor (*), more formally stated: $G/N_{k-1}\cong A_5$, because otherwise if there were a chief series with $G/N_{k-1}$ not being isomorphic to $A_5$, then the action of $N_{k-1}$ on $V$ would be a smaller counterexample, a contradiction.\\

We now prove that $V_N$ (which denotes $V$ viewed as an $N$-module) is homogeneous. By Clifford we can write $V_N$ as a direct sum of its homogeneous components. Suppose that there are $l$ of them and that $l>1$. Then $G$ permutes them transitively, and if $L$ is the kernel of this permutation action, then by (*), $G/L$ has a chief factor isomorphic to $A_5$. In particular, $G/L$ contains an element of order 3 which will have an orbit of size 3 in this permutation action, thus leading it to have a nontrivial fixed point on $V$, contradicting the hypothesis. Hence $l=1$ which means that $V_N$ is homogeneous.\\

Now write $G_0=G/C_G(V)$ and $N_0=N/C_N(V)=N/C_G(V)$. 
We want to show that $N_0$ is a 2-group. For this, recall that 3 and 5 do not divide $|N_0|$ and suppose that $r>5$ is a prime dividing $|N_0|$. Then there exist normal subgroups $S$, $T$ of $G_0$ such that $S/T$ is a chief factor of $G_0$ and $r$ does not divide $|G_0/S|$. Then by Schur-Zassenhaus, $S/T$ is a complemented (or non-Frattini) chief factor, and so there exists a subgroup $K_0$ of $G_0$ such that $T\leq K_0$, $K_0S=G_0$, and $(K_0/T)\cap (S/T)=1$. So if $K$ is the inverse image of $K_0$ in $G$, then the action of $K$ on $V$ is a smaller counterexample, contradicting our minimal choice. Hence $r$ does not exist and $N_0$ is a 2-group.\\

Next we show that $\phi(N_0)=1$ where $\phi(\cdot)$ denotes the Frattini subgroup. Assume that $\phi(N_0)>1$. Clearly the Fitting subgroup $F(G_0)$ of $G_0$ equals $N_0$, and $F(G_0)/\phi(G_0)$ has a complement in $G_0/\phi(G_0)$. So let $H_0\leq G_0$ such that $\phi(N_0)\leq H_0$ and $H_0/\phi(G_0)$ is such a complement of $F(G_0)/\phi(G_0)$. If $H$ is the preimage of $H_0$ in $G$, then the action of $H$ on $V$ is a smaller counterexample, contradicting our minimal choice. This establishes that indeed $\phi(N_0)=1$.\\

Thus $N_0$ is an elementary abelian 2-group. Since $V_N$ is homogeneous (as seen above), then also $V_{N_0}$ is homogeneous, and since $N_0$ is elementary abelian and acts faithfully on $V$, we see that altogether $N_0$ must be cyclic of order 2. This forces $N_0$ to be the center of $G_0$, and thus $G_0$ is a central extension of $A_5$, which by  our hypothesis and 
Lemma \ref{plms lemma}, $G_0$ must be non-split. It is well-known that the only such extension is SL(2,5). Therefore the action of $G_0$ on $V$ is a faithful irreducible coprime action. By \cite[Lemma 10]{robinson-thompson} we may assume that $V$ is an absolutely reducible $G_0$-module, and since the action is coprime, it is well-known that $V$ can be viewed as a module over the complex numbers. Hence we need to study the complex irreducible faithful representations of $\SL(2,5)$ and check that if the elements of order 3 act frobeniusly on the underlying module, then the elements of order 5 also act frobeniusly on that module.\\

Now SL(2,5) has four faithful representations. Two of them are of dimension 2, and in these SL(2, 5) acts frobeniusly on the 
underlying module. In particular, all 3-elements and 5-elements act frobeniusly, which is in line with our assertion (or rather, which
leads to a contradiction to $G$, $V$ being a counterexample, as desired). Also, the faithful representation of degree 6 of $\SL(2,5)$ is
just induced from a subgroup of order 20, and thus, in particular, the elements of order 3 have nontrivial fixed points on $V$,
contradicting the hypothesis. The remaining faithful representation is of degree 4 is in some ways the hardest to understand.
Anyway, we used GAP \cite{gap} to verify that here again, the elements of order 3 have nontrivial fixed points on the corresponding module,
again contradicting the hypothesis. This final contradiction shows that the minimal counterexample does not exist and completes
the proof of the lemma. \end{proof}


\end{document}